\magnification 1200
\let\nd\noindent 
\font\tenmsb=msbm10   
\font\sevenmsb=msbm7
\font\fivemsb=msbm5
\newfam\msbfam
\textfont\msbfam=\tenmsb
\scriptfont\msbfam=\sevenmsb
\scriptscriptfont\msbfam=\fivemsb
\def\Bbb#1{\fam\msbfam\relax#1}
\def\Z{{\Bbb Z}}
\def\Q{{\Bbb Q}}
\def\R{{\Bbb R}}
\def\F{{\Bbb C}}
\def\A{{\cal A}}
\def\B{{\cal B}}
\def\C{{\cal C}}

\def\G{{\cal G}}

\def\Cech{{\v Cech }}
\def\im{{\rm Im}}
\def\star{{\rm star}}
\def\ke{{\rm Ker}}

\def\qed{\hbox{\hskip 6pt\vrule width6pt height7pt depth1pt \hskip1pt}}

\def\vs#1 {\vskip#1truein}
\def\hs#1 {\hskip#1truein}
  \hsize=6truein        \hoffset=.25truein 
  \vsize=8.8truein  
  \pageno=1     \baselineskip=12pt
  \parskip=0 pt         \parindent=20pt
  \overfullrule=0pt     \lineskip=0pt   \lineskiplimit=0pt
  \hbadness=10000 \vbadness=10000 

%
%
\newdimen\FigSize	\FigSize=.9\hsize 
%
\newskip\abovefigskip	\newskip\belowfigskip
\gdef\epsfig#1;#2;{\par\vskip\abovefigskip\penalty -500
   {\everypar={}\epsfxsize=#1\noindent
    \centerline{\epsfbox{#2}}}%
    \vskip\belowfigskip}%
%
\newskip\figtitleskip
\gdef\tepsfig#1;#2;#3{\par\vskip\abovefigskip\penalty -500
   {\everypar={}\epsfxsize=#1\noindent
    \vbox
      {\centerline{\epsfbox{#2}}\vskip\figtitleskip
       \centerline{\figtitlefont#3}}}%
    \vskip\belowfigskip}%
%
\newcount\FigNr	\global\FigNr=0
\gdef\nepsfig#1;#2;#3{\global\advance\FigNr by 1
   \tepsfig#1;#2;{Figure\space\the\FigNr.\space#3}}%
%
%
%
\gdef\ipsfig#1;#2;{
   \midinsert{\everypar={}\epsfxsize=#1\noindent
	      \centerline{\epsfbox{#2}}}%
   \endinsert}%
%
\gdef\tipsfig#1;#2;#3{\midinsert
   {\everypar={}\epsfxsize=#1\noindent
    \vbox{\centerline{\epsfbox{#2}}%
          \vskip\figtitleskip
          \centerline{\figtitlefont#3}}}\endinsert}%
%
\gdef\nipsfig#1;#2;#3{\global\advance\FigNr by1%
  \tipsfig#1;#2;{Figure\space\the\FigNr.\space#3}}%
\newread\epsffilein    
\newif\ifepsffileok    
\newif\ifepsfbbfound   
\newif\ifepsfverbose   
\newdimen\epsfxsize    
\newdimen\epsfysize    
\newdimen\epsftsize    
\newdimen\epsfrsize    
\newdimen\epsftmp      
\newdimen\pspoints     
\pspoints=1bp          
\epsfxsize=0pt         
\epsfysize=0pt         
\def\epsfbox#1{\global\def\epsfllx{72}\global\def\epsflly{72}%
   \global\def\epsfurx{540}\global\def\epsfury{720}%
   \def\lbracket{[}\def\testit{#1}\ifx\testit\lbracket
   \let\next=\epsfgetlitbb\else\let\next=\epsfnormal\fi\next{#1}}%
\def\epsfgetlitbb#1#2 #3 #4 #5]#6{\epsfgrab #2 #3 #4 #5 .\\%
   \epsfsetgraph{#6}}%
\def\epsfnormal#1{\epsfgetbb{#1}\epsfsetgraph{#1}}%
\def\epsfgetbb#1{%
%
%
\openin\epsffilein=#1
\ifeof\epsffilein\errmessage{I couldn't open #1, will ignore it}\else
%
%
   {\epsffileoktrue \chardef\other=12
    \def\do##1{\catcode`##1=\other}\dospecials \catcode`\ =10
    \loop
       \read\epsffilein to \epsffileline
       \ifeof\epsffilein\epsffileokfalse\else
%
%
          \expandafter\epsfaux\epsffileline:. \\%
       \fi
   \ifepsffileok\repeat
   \ifepsfbbfound\else
    \ifepsfverbose\message{No bounding box comment in #1; using defaults}\fi\fi
   }\closein\epsffilein\fi}%
%
%
\def\epsfsetgraph#1{%
   \epsfrsize=\epsfury\pspoints
   \advance\epsfrsize by-\epsflly\pspoints
   \epsftsize=\epsfurx\pspoints
   \advance\epsftsize by-\epsfllx\pspoints
%
%
   \epsfxsize\epsfsize\epsftsize\epsfrsize
   \ifnum\epsfxsize=0 \ifnum\epsfysize=0
      \epsfxsize=\epsftsize \epsfysize=\epsfrsize
%
%
     \else\epsftmp=\epsftsize \divide\epsftmp\epsfrsize
       \epsfxsize=\epsfysize \multiply\epsfxsize\epsftmp
       \multiply\epsftmp\epsfrsize \advance\epsftsize-\epsftmp
       \epsftmp=\epsfysize
       \loop \advance\epsftsize\epsftsize \divide\epsftmp 2
       \ifnum\epsftmp>0
          \ifnum\epsftsize<\epsfrsize\else
             \advance\epsftsize-\epsfrsize \advance\epsfxsize\epsftmp \fi
       \repeat
     \fi
   \else\epsftmp=\epsfrsize \divide\epsftmp\epsftsize
     \epsfysize=\epsfxsize \multiply\epsfysize\epsftmp   
     \multiply\epsftmp\epsftsize \advance\epsfrsize-\epsftmp
     \epsftmp=\epsfxsize
     \loop \advance\epsfrsize\epsfrsize \divide\epsftmp 2
     \ifnum\epsftmp>0
        \ifnum\epsfrsize<\epsftsize\else
           \advance\epsfrsize-\epsftsize \advance\epsfysize\epsftmp \fi
     \repeat     
   \fi
%
%
   \ifepsfverbose\message{#1: width=\the\epsfxsize, height=\the\epsfysize}\fi
   \epsftmp=10\epsfxsize \divide\epsftmp\pspoints
   \vbox to\epsfysize{\vfil\hbox to\epsfxsize{%
      \includegraphics{#1}%
      \hfil}}%
\epsfxsize=0pt\epsfysize=0pt}%
%
%
{\catcode`\%=12 \global\let\epsfpercent=
%
%
\long\def\epsfaux#1#2:#3\\{\ifx#1\epsfpercent
   \def\testit{#2}\ifx\testit\epsfbblit
      \epsfgrab #3 . . . \\%
      \epsffileokfalse
      \global\epsfbbfoundtrue
   \fi\else\ifx#1\par\else\epsffileokfalse\fi\fi}%
%
%
\def\epsfgrab #1 #2 #3 #4 #5\\{%
   \global\def\epsfllx{#1}\ifx\epsfllx\empty
      \epsfgrab #2 #3 #4 #5 .\\\else
   \global\def\epsflly{#2}%
   \global\def\epsfurx{#3}\global\def\epsfury{#4}\fi}%
%
%
\def\epsfsize#1#2{\epsfxsize}%
%
%

\epsfverbosetrue			
\abovefigskip=\baselineskip		
\belowfigskip=\baselineskip		
\global\let\figtitlefont\bf		
\global\figtitleskip=.5\baselineskip	

\pageno=0

\footline{\ifnum\pageno=0\hss\else\hss\tenrm\folio\hss\fi}
\hbox{}
\vskip 1truein\centerline{{\bf A HOMEOMORPHISM INVARIANT FOR SUBSTITUTION }}
\vskip .1truein\centerline{{\bf TILING SPACES}}
\vskip .2truein\centerline{by}
\vskip .2truein
\centerline{Nicholas Ormes${}^{1,2}$,
\ \ {Charles Radin${}^{1}$ 
\footnote{*}{Research supported in part by Texas ARP Grants 003658-152
and 003658-158\hfil}}\ \ and {Lorenzo Sadun${}^{1,3}$
\footnote{**}{Research supported in part by Texas ARP Grants 003658-152
and 003658-158 
\hfill\break \indent and NSF Grant DMS-9626698\hfil}}}
\vskip .5truein\centerline{\vbox{
${}^1$\ \ Mathematics Department, University of Texas at Austin
\vskip.1truein
${}^2$\ \ Mathematics Department, University of Connecticut
\vskip.1truein
${}^{3}$\ \ Physics Department, Technion -- Israel Institute of
Technology}}
\vs.5
\centerline{{\bf Abstract}}
\vs.1 \nd
We derive a homeomorphism invariant for those tiling spaces which are made
by rather general substitution rules on polygonal tiles, including
those tilings, like the pinwheel, which contain tiles in infinitely
many orientations. The invariant is a quotient of \Cech cohomology,
is easily computed directly from the substitution rule, and
distinguishes many examples, including most pinwheel-like tiling
spaces. We also introduce a module structure on cohomology which is
very convenient as well as of intuitive value.
\vs.8
\centerline{July 2000}
\vs.2
\centerline{Subject Classification:\ \ 37B50, 52C23, 52C20}
\vfill\eject 

\nd {\bf 1. Introduction.}

In this paper we study the topology of substitution tiling
spaces. Specifically, we are interested in ways in which 
the geometry of the substitution is reflected in the topology. 

We follow a history of work where one studies the topology of spaces
created with an underlying dynamical system. An early reference to the
circle of ideas is Parry and Tuncel [PaT; chap. 4]. They show that the first
cohomology group of the standard suspension space $\widetilde{X}$ for a
$\Bbb Z$-subshift $(X,T)$ is a useful invariant for topological
conjugacy of $(X,T)$.  Here, $\widetilde{X} = (X \times [0,1]) /\sim$, where
$(x,1) \sim (Tx,0)$.

Later, Herman, Putnam and Skau [HPS] associated a C$^\ast$-algebra
with general minimal $\Bbb Z$-subshifts via the crossed product
construction, in analogy with the way von Neumann algebras had been
associated previously with measureable dynamics. For these systems the
$K_0$-group of the C$^\ast$-algebra coincides with the first
cohomology of the suspension space [HPS]. 
Giordano, Putnam and Skau [GPS] used the (ordered) $K$-groups of
the C$^\ast$-algebra to characterize orbit equivalence for minimal
$\Bbb Z$-subshifts.

In a similar vein Connes [Con] used Penrose tilings to motivate 
ideas on noncommutative topology. Specifically, he considered the
quotient $X$ of the space of Penrose tilings by the equivalence
relation of translation. Using the hierarchical structure of the
Penrose tilings there is a natural homeomorphism between $X$ and a
quotient $Y$ of sequences modulo another equivalence relation. As is
common for equivalence relations coming from dynamics, the quotient
topology of $Y$ is trivial.  And so Connes considers a
(noncommutative) C$^\ast$-algebra associated with the quotient
$Y$. The $K$-theory of this C$^\ast$-algebra is his substitute for the
(trivial) quotient topology.

This was pushed one step further by Kellendonk and Putnam [Ke1-4, KeP,
AnP] in which examples of tiling dynamical systems such as the Penrose
system are considered and their (unordered) $K$-theory is worked out,
obtaining for instance cohomology groups for the spaces of tilings.

We consider substitution tiling systems in $\R^d$, such as the Penrose
tilings in $\R^2$, as dynamical systems with $\R^d$ (translation)
actions. These are natural generalizations of the suspensions of
$\Z^d$ (substitution) subshifts. 
(We will restrict attention to $d=2$; while most of
our methods carry over to higher dimensions, there are complications which
we wish to avoid.)  As noted above, suspensions are
sometimes introduced to get invariants to distinguish between
subshifts, insofar as the spaces on which the subshifts act are all
the same --- Cantor sets.  For tilings the natural objects already
have $\R^d$ actions and different tiling spaces need not be
homeomorphic.  In particular, the cohomology groups of such spaces are
invariants of topological conjugacy. For 1-dimensional tiling systems,
Barge and Diamond [BD] worked out complete homeomorphism invariants.

Although in principle it should be possible to work out the cohomology
groups for tiling spaces, in practice it can be difficult since the
number of different tiles in interesting models (especially after
collaring --- see below), which gives the number of generators for the
chain groups, is usually large, even infinite for examples such as
pinwheel tilings. New in this work is analysis of two aspects of the
cohomology of substitution tiling spaces that are both easy to compute
and informative about the topology of the spaces.

The first aspect is an order structure for the top dimensional \Cech
cohomology. For some history in a ``noncommutative'' setting see
[GPS], but the idea is much older, appearing for instance in
[SuW]. Fundamental to this paper as well is a representation of our
tiling spaces as inverse limits, as discussed in [AnP], and in particular
we use this to define a positive cone in cohomology.
More specifically, we construct a homeomorphism
from the top \Cech cohomology group to the real numbers, and define the positive cone
to be the preimage of the non-negative reals. The image of the map is
closely related to the additive group of $\Z[1/\lambda]$, where
$\lambda$ is the area stretching factor of the substitution.  We
show that the positive cone is a homeomorphism invariant, and that one
can therefore extract invariant topological information from the
easily computable stretching factor $\lambda$.  

We prove results of this form for two classes of tiling spaces.  For
``fixed-orientation'' tiling spaces (denoted $X_x$ below), we need the
additional assumption that tiles only appear in a finite number of
distinct orientations in each tiling. For ``all-orientation''
tiling spaces (denoted $X_\phi$ below), this assumption is not
needed. Our main result is
\vs.1
\nd {\bf Theorem 1.} {\it Let $X$ and $Y$ be substitution tiling systems whose
substitutions have area stretching factors $\lambda$ and $\lambda'$,
respectively.  Suppose that either i) $X$ and $Y$ have finite relative
orientation groups and the fixed-orientation spaces $X_x$ and $Y_y$
are homeomorphic, or ii) the all-orientation spaces $X_\phi$ and
$Y_\phi$ are homeomorphic. Then

\item{(1)} If $\lambda$ is an integer so is $\lambda'$, and 
$\Z[1/\lambda] =\Z[1/\lambda']$ as subsets of $\R$.  In particular,
$\lambda$ and $\lambda'$ have the same prime factors, although not
necessarily with the same multiplicities.

\item{(2)} If $\lambda$ is not an integer then $\lambda$ and $\lambda'$ 
are irrational, and $\Q[\lambda] = \Q[\lambda']$ 
as subsets of $\R$.}
 
\vs.1
(We will see in examples below that this distinguishes topologically
between tiling spaces, such as that of the pinwheel and the
(2,3)-pinwheel, which had not previously been known to be
distinguishable even in the stronger sense of topological conjugacy, that
is, including the natural $\R^2$ action of translations on the spaces.)

To prove this theorem we show that the order structure, although
defined in terms of the (not necessarily invariant) inverse limit
structure, is in fact invariant.  The relations between $\lambda$ and
$\lambda'$ then follow by algebraic arguments.

To describe the second aspect of our analysis we must discuss
the unusual rotational properties of substitution tilings such as those
of Penrose. No Penrose tiling is invariant under rotation by $2\pi/10$
about any point, but every translation invariant Borel probability
measure on the space of Penrose tilings ($X_x$ or $X_\phi$)
is invariant under rotation by
$2\pi/10$. We say $\Z_{10}$ is the ``relative orientation group'' of
the Penrose tilings.  This symmetry has been the source for most of
the explosion of work on such tilings, in particular their use in
modelling quasicrystals [Ra2]. (We have already investigated the use of
rotational symmetry as a conjugacy invariant [RaS].)

As we will see below, each element of the cyclic rotation group
$\Z_{10}$ gives an automorphism of the space $X_x$ of Penrose tilings.
As a result, the group $\Z_{10}$, and by extension its group ring
$\Z[t]/(t^{10}-1)$, acts naturally on the (\Cech) cohomology groups of
$X_x$.  This gives $H^*(X_x)$ the structure of a module over the group
ring, and this module structure is quite useful.  In particular, for
any realization of $X_x$ as the inverse limit of simplicial complexes,
as in [AnP], the (co)chain groups are themselves modules, and the
boundary and substitution maps are equivariant.  This greatly
simplifies the calculations.  For Penrose tiles for example [AnP],
there are 4 types of tiles, each appearing in 10 different
orientations.  Without using the rotational symmetry the action of
substitution on 2-chains would be given by a $40 \times 40$ matrix.
Using the rotational symmetry we express it instead as a $4 \times 4$
matrix with entries in the group ring of $\Z_{10}$, which can be
analyzed one irreducible representation at a time. The substitutions
on 0- and 1-chains, and the boundary maps, are handled similarly.  The
result is a streamlined calculation that yields information beyond
that obtained by older methods.  Anderson and Putnam [AnP] computed
$H^0(X_x)=\Z$, $H^1(X_x)=\Z^5$ and $H^2(X_x)=\Z^8$.  We note that, as
modules, a finite-index additive subgroup of $H^1(X_x)$ equals
$\Z[t]/(t-1) \oplus \Z[t]/(t^4-t^3+t^2-t+1)$ while a finite index
subgroup of $H^2(X_x)$ equals $(\Z[t]/(t-1))^2 \oplus (\Z[t]/(t+1))^2
\oplus \Z[t]/(t^4-t^3+t^2-t+1)$, where rotation by $2\pi/10$ is
multiplication by $t$.

\vs.2
\noindent{\bf 2. Tiling spaces and their associated inverse limits.} 

Before defining substitution tiling spaces 
in general we present some examples whose topological properties we will
later consider.

A ``chair'' tiling of the plane, Fig.~1, can be made as follows. Consider the
L-shaped tile of Fig.~2. Divide this tile (also called a ``tile of
level 0'') into four pieces as in Fig.~3 and rescale by a linear
factor of 2 so that each piece is the same size as the original.  This
yields a collection of 4 tiles that we call a ``tile of level 1''.
Subdividing each of these tiles and rescaling gives a collection of 16
tiles that we call a tile of level 2.  Repeating the process $n$ times
gives a collection of $4^n$ tiles --- a tile of level $n$.  A ``chair''
tiling is a tiling of the plane with the property that every finite
subcollection of tiles is congruent to a subset of a tile of some
level.  A chair tiling has only one type of tile, appearing in 4
different orientations. (For chair tilings, and in fact very generally
for tilings made by such substitution rules, there are ``matching rules''
which provide a different method of construction [Goo, Ra2].)

Somewhat more complicated are the ``Penrose'' tilings, Fig.~4, which have 4 different
types of triangular tiles, each appearing in 10 different orientations.  The 4 tiles, and
the rule for subdividing them, are shown in Fig.~5. (These tiles are not the
familiar kites and darts.  However, the space of substitution tilings developed 
from these 4 triangles is homeomorphic to the space of Penrose kite and dart tilings).
Here the linear stretching factor is the golden mean $\tau=(1+\sqrt{5})/2$, and the
area 
\vfill \eject
\hbox{}\vs1
\vbox{\epsfig 1.2\hsize; 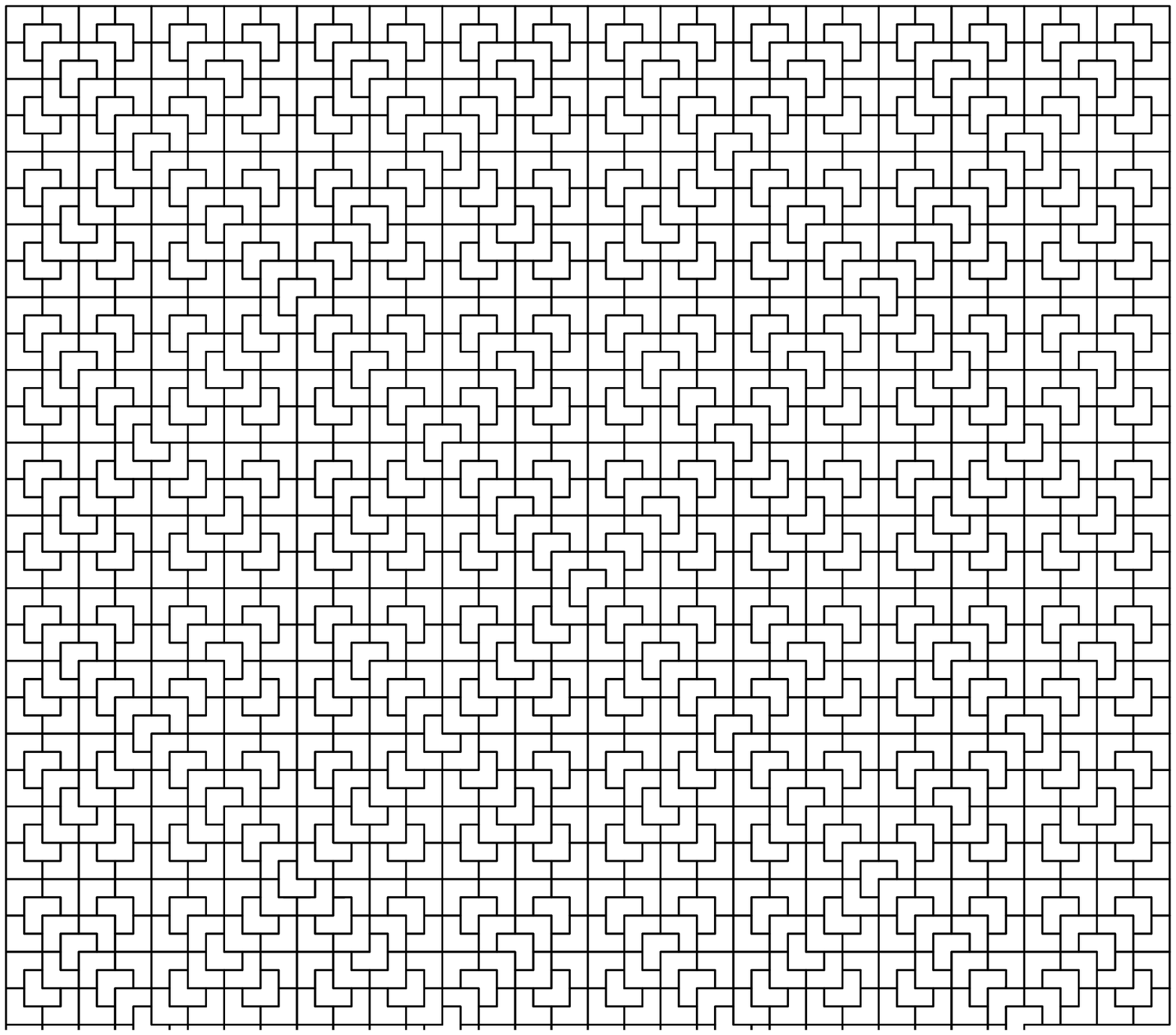 ;}
\vs.5
\centerline{Figure 1. A chair tiling}
\vfill \eject
\hbox{}
\vbox{\epsfig .3\hsize; 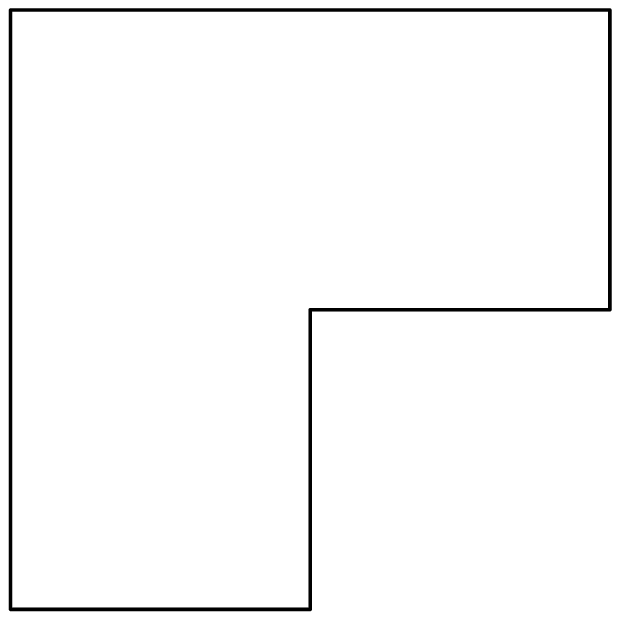;}
\vs.3
\centerline{Figure 2. The chair tile}
\vs.8
\hs1.5 \vbox{\epsfig 1\hsize; 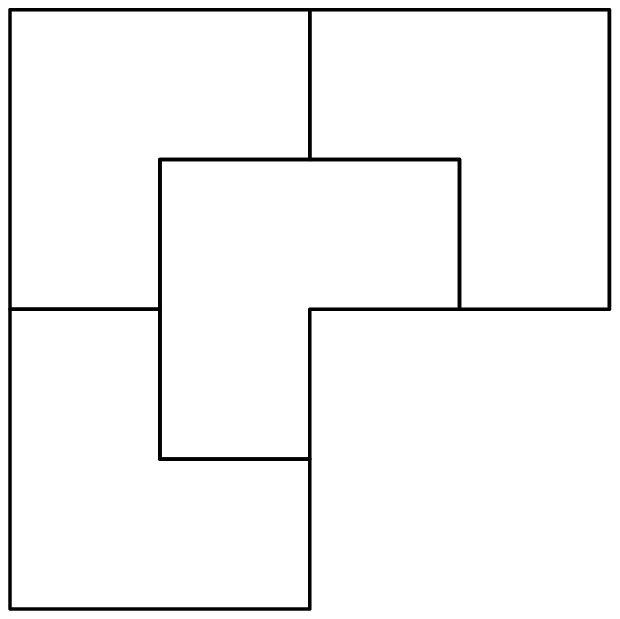 ;}
\vs-1.6
\centerline{Figure 3. The chair substitution}
\vfill \eject
\hbox{}
\vbox{\epsfig .8\hsize; 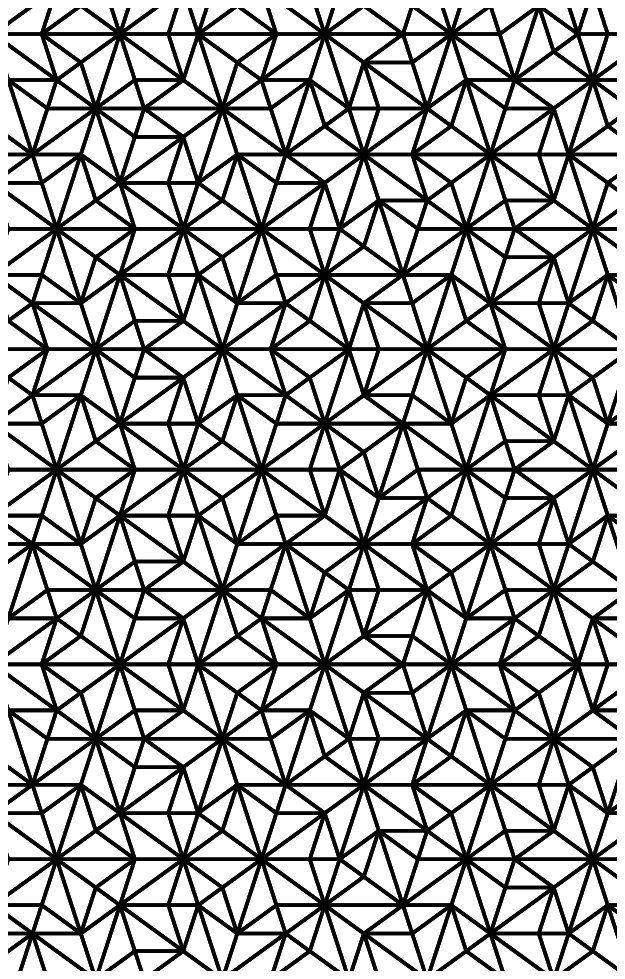;}
\vs.3
\centerline{Figure 4. A Penrose tiling}
\vfill\eject
\hbox{}
\vbox{\epsfig 1\hsize; 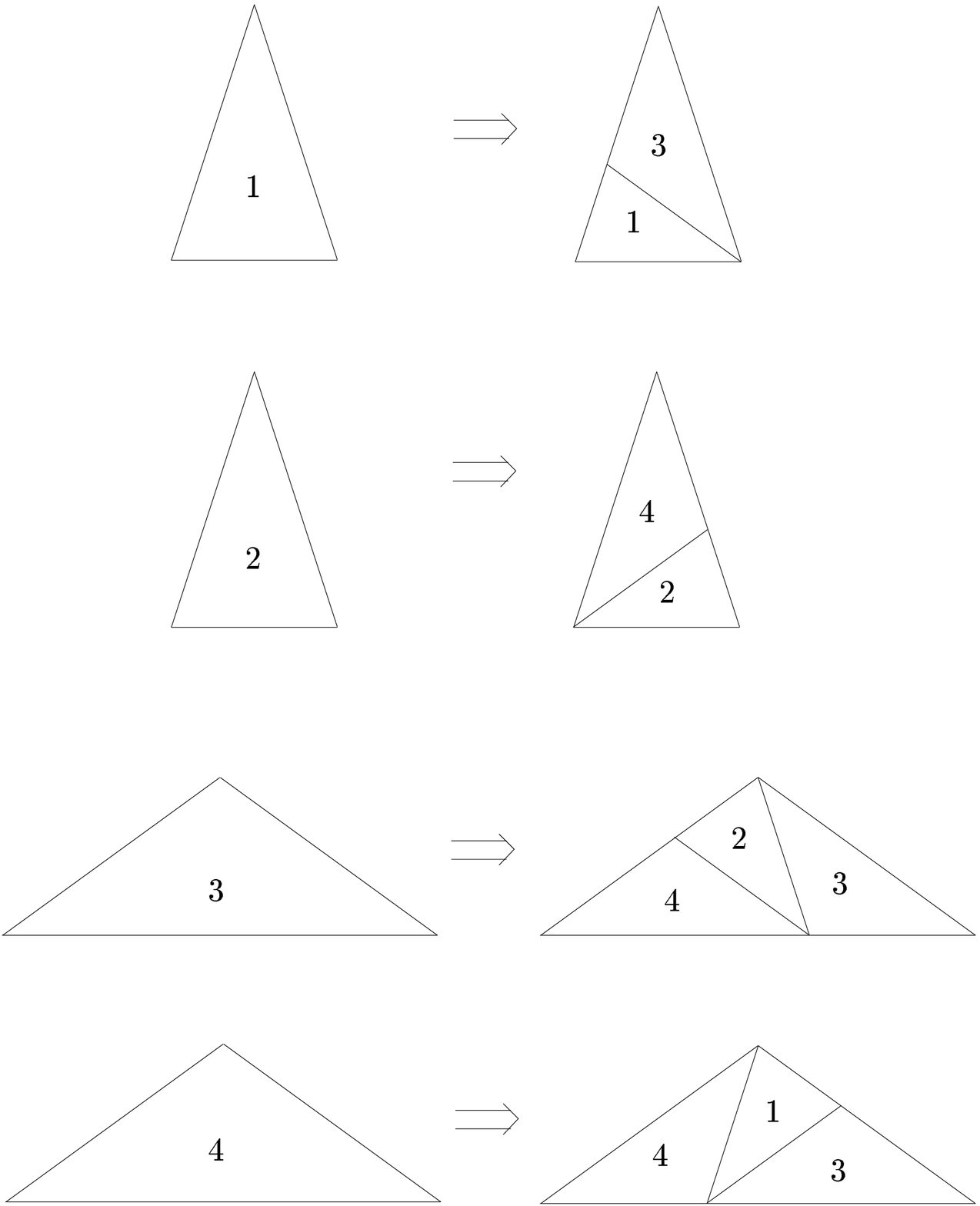;}
\vs.3
\centerline{Figure 5. The Penrose substitution}
\vfill\eject
\hbox{}\vs-1
\hs.15 \epsfig 1.1\hsize; 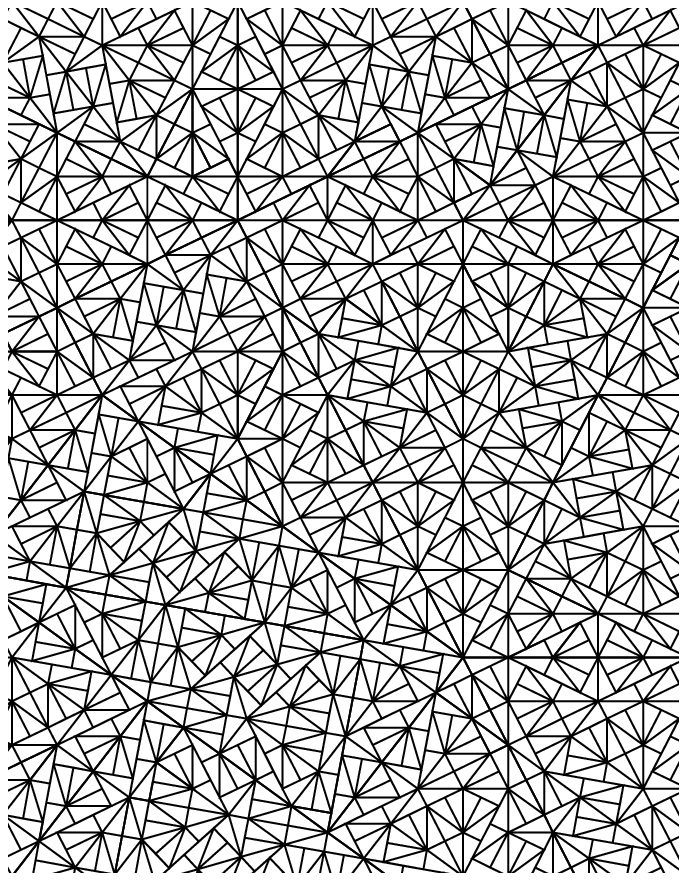;
\vs0 \centerline{Figure 6. A pinwheel tiling}
\vfill \eject
\hbox{}\vs1
\vbox{\epsfig 1\hsize; 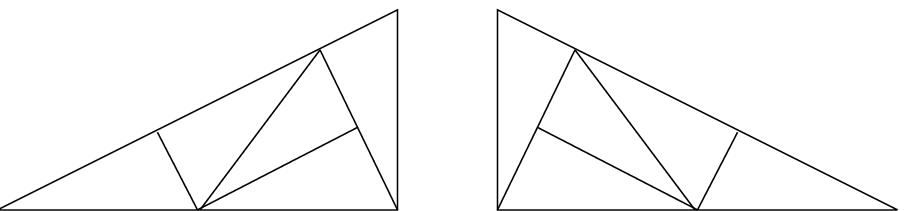 ;}
\vs.3
\centerline{Figure 7. The substitution for pinwheel tilings}
\vs1
\vbox{\epsfig 1\hsize; 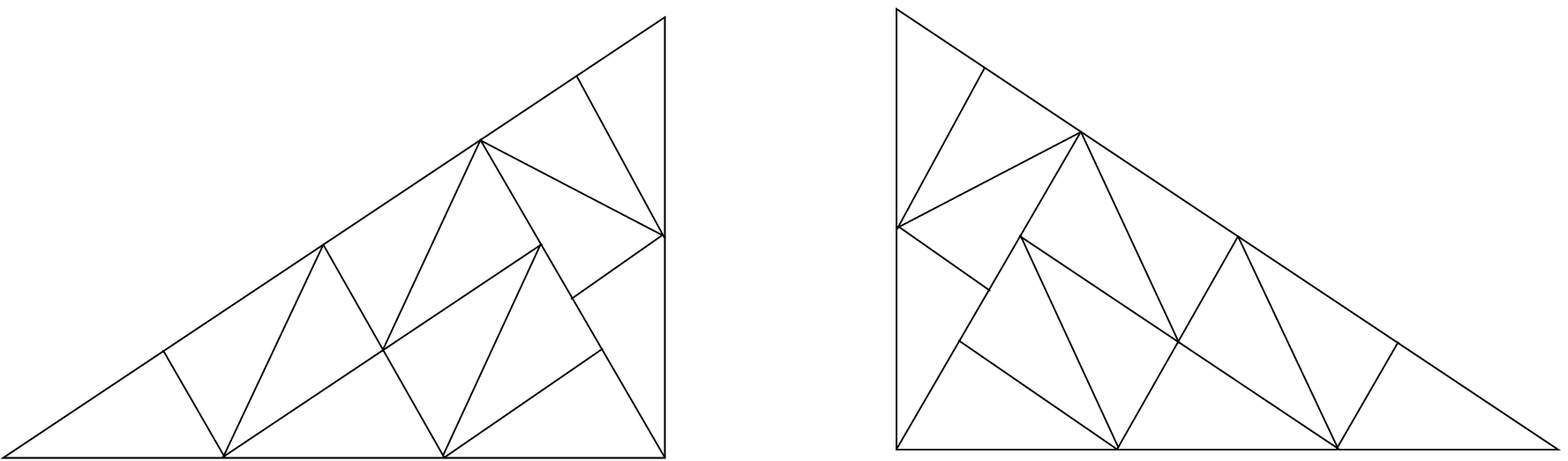 ;}
\vs.3
\centerline{Figure 8. The substitution for (2,3) pinwheel tilings}
\vfill \eject
\nd
stretching factor is $\tau^2=(3+\sqrt{5})/2$.

A ``pinwheel'' tiling of the plane, Fig.~6 [Ra1], has tiles that appear in an
infinite number of orientations.  The two basic tiles, a
1-2-$\sqrt{5}$ right triangle and its mirror image, are shown in
Fig.~7, with their substitution rule.  

Notice that at the center of a tile of level 1 there is a tile of
level 0 similar to the level 1 tile but rotated by an angle $\alpha =
\tan^{-1}(1/2)$. Similarly, the center tile of a tile of level $n$ is
rotated by $n\alpha$ relative to the tile.  Since $\alpha$ is an
irrational multiple of $\pi$, we see, using the fact that within a
tile of level 2 there is a tile of level 0 similar to the level 2
tile, that this rotation never stops, and each tiling contains tiles
in an infinite number of distinct orientations.

Finally, we consider the ``(2,3)-pinwheel'' tilings defined by the
substitution of Fig.~8, whose tiles are
2-3-$\sqrt{13}$ right triangles. (Such variants of the pinwheel 
are easily constructed
for any integral legs $m<n$.)  Like the ordinary pinwheel, variant
pinwheel tilings also necessarily have tiles in an infinite number of
distinct orientations. In fact the relative orientation groups for
all pinwheel tilings are algebraically isomorphic. Theorem 1 shows
that the tiling spaces for the pinwheel and (2,3)-pinwheel are not 
homeomorphic.

In all the above cases it is easy to construct explicit examples of
tilings.  Pick a tile to include the origin of the plane.  Embed this
tile in a tile of level 1 (there are several ways to do this).  Embed
that tile of level 1 in a tile of level 2, embed that in a tile of
level 3, and so on.  The union of these tiles will cover an infinite
region, typically -- though not necessarily -- the entire plane.
 
We now give a general definition of substitution tiling systems.
Let $\A$ be a nonempty finite collection of polygons in the plane.
Let $X(\A)$ be the set of all tilings
of the plane by congruent copies, which we call tiles, of
the elements of (the ``alphabet'') $\A$. We label the 
``types'' of tiles by the elements of $\A$.  We endow $ X(\A)$ with the metric
$$
d(x, x')\equiv \sup_{n}{1\over n}m_{H}[B_{n}(\partial x),B_{n}(\partial
x')], \eqno 1)
$$
where ${B}_{n}(\partial x)$ denotes the intersection of two sets: the
closed disk ${B}_n$ of radius $n$ centered at the origin of the
plane and the union $\partial x$ of the boundaries 
of all tiles in $x$. The Hausdorff metric $m_{H}$ is defined as
follows. Given two compact subsets $P$ and $Q$ of $\R^2$, $m_{H}[P,Q]
= \max \{ {\tilde d} (P,Q), {\tilde d} (P,Q)\}$, where
$${\tilde d} (P,Q) =  \sup_{p \in P} \inf_{q \in Q} ||p - q||, \eqno 2)$$
with $||w||$ denoting the usual Euclidean norm
of $w$.  

Under this metric two tilings are  close if they agree, up to a small
Euclidean motion, on a large  disk around the origin.  The converse is
also true for  tiling systems whose tiles meet full  edge to full edge
(as we require) -- closeness  implies agreement, up to small Euclidean
motion, on  a large disk around  the origin.  Although  the metric $d$
depends on the location of the  origin, the topology induced by $d$ is
translation invariant.  A sequence  of tilings converges in the metric
$d$ if and  only if its restriction to every  compact subset of $\R^2$
converges in $m_H$.  It is not  hard to show [RaW] that $X(\A)$ 
is compact and that the
natural action  of the connected  Euclidean group $\G_E$ on  $ X(\A)$,
$(g,x)\in   \G_E\times   X(\A)\longrightarrow  g[x]\in   X(\A)$,   is
continuous.

A ``substitution tiling space'' is a closed subset $X_\phi \subset
X(\A)$ satisfying some additional conditions. 
To understand these conditions we first need the notion of
``patches''.  A patch is a (finite or infinite) subset of an element
$x\in X(\A)$; the set of all patches for a given alphabet will be
denoted by $W$. Next we need, as for the above examples, an auxiliary
``substitution function'' $\phi$, a map from $W$ to $W$, with the
following properties:
\vs.1
\item{(1)} There is some constant $c(\phi)>1$ such 
that, for any $g\in \G_E$ and $x\in X$, 
$\phi(g[x])=\phi(g)[x]$, where 
$\phi(g)$ is the conjugate of $g$ by the similarity of
Euclidean space consisting of
stretching about the origin by $c(\phi)$.
\item{(2)} For each tile $T\in\A$ and for each 
$n \ge 1$, the union of the tiles in $\phi^n T$ is
congruent to $[c(\phi)]^n T$, and these
tiles meet full edge to full edge.

\item{(3)}  For each tile $T\in\A$, $\phi T$ contains at
least one tile of each type.

\item{(4)} For each tile $T\in\A$ there is $n_T \ge 1$
such that $\phi^{n_T} a$ contains a tile of
the same type as $T$ and 
parallel to it.
\item{(5)} No tile $T\in\A$ has a nontrivial rotational symmetry.
\vs.1 

Condition (2) is significant. It is satisfied by the pinwheel tilings
only if we add additional vertices at midpoints of the legs of length
2, creating boundaries of 4 edges. (A similar alteration is needed for the
chair tilings.) The tile of level $n$, $\phi^n T$, will be
said to be of ``type'' $T$.

Condition (5), by contrast, is technical, and does not significantly limit 
the scope of this
work.  If a tile does have $n$-fold rotational symmetry,
we can recover condition (5) by breaking the tile into $n$ congruent but asymmetric
pieces, in a manner consistent with the subdivision rules.
\vs.1
\nd {\bf Definition 1.} For a given alphabet $\A$ of polygons and
substitution function $\phi$  the ``substitution tiling space''
is the compact subspace $X_\phi\subset X(\A)$
of those tilings $x$ such that every finite
subpatch of $x$ is congruent to a subpatch of $\phi^n(T)$ for some $n>0$ and
$T\in \A$. We assume $x$ can be decomposed in one and only
one way by tiles of level $n$, for any fixed $n\ge 0$ [Sol].

\vs.1
The above definition gives spaces of tilings in which the tiles can appear
in arbitrary orientations. Although any {\it fixed} chair tiling, 
for example, has tiles in only 
four orientations, the space of {\it all} chair tilings also contains rotated 
versions of that tiling, and so contains tilings in which chairs appear in any orientation.

In many cases, especially when working with tilings whose tiles have
only a finite number of orientations per tiling, it is convenient (and
customary) to allow tiles only a minimal set of orientations.  Instead
of using the metric $d$, we define a metric $d'$ such that two tilings
are close if (and only if) they agree on a large disk around the
origin {\it up to a small translation}, rather than up to a small
Euclidean motion. Specifically, 
$$d'(x,x')=\inf\{{R+S\over 1+R+S}:B_{1+1/S}(g[\partial x])=B_{1+1/S}(h[\partial
x'])\} \eqno 2')$$
\nd where $g,h\in\R^2$, 
$R>0,\ S>0,\ |g|<R,\ |h|<R$ and $d'(x,x')\equiv 1$ if the defining set is empty.
We pick any one tiling $ x \in X_\phi$, and define
$X_{ x}\subset X_\phi$ to be $\{x'\in X_\phi: d'(x,x')<1\}$.
For tiling spaces
with only a finite number of orientations per tiling, $X_{ x}$
consists of those tilings in which tiles only appear parallel to tiles
in $ x$. We call $X_{ x}$ a ``fixed-orientation'' tiling space.  

Finally, we consider the quotient of $X_\phi$ by rotations about the
origin and denote this space $X_0$.  For any fixed tiling $ x$, $X_0$
is also the quotient of $X_{ x}$ by the ``relative orientation group''
$G$ of the tiling $ x$, defined as follows: $G$ is generated by the
rotations that take a tile in $ x$ and turn it parallel to another
tile of the same type in $ x$.  In [RaS] it was shown that this group
is the same for all tilings $ x \in X_\phi$. The
topologies of these three spaces are closely related. 
The space $X_\phi$ is a circle bundle over
$X_0$ with some singular fibers; these fibers correspond to tilings
with a discrete rotational symmetry about the origin. 

Next we show how to give a substitution tiling space the structure of
an inverse limit.  We begin by doing this first for a fixed-orientation 
space, and we begin in particular with the example of the chair.  Fix a
chair tiling $x$, and in it a tile $T$.  Each of the other tiles in
$x$ are rotated with respect to $T$ by an element of the cyclic
rotation group $G_{ch} = \Z_4$. We construct an (uncollared) chair
complex $\Sigma^{uc}$ as follows.  We start with the topological
disjoint sum $\Sigma'$ of the elements of $\{g[T]\,:\, g\in G_{ch}\}$.
We then identify those edges in $\Sigma'$ which ``meet'' somewhere in
$x$, defining $\Sigma^{uc}$.  We consider a countable number of copies
$\Sigma^{uc}_j$ of $\Sigma^{uc}$, indexed by the non-negative
integers. We think of the tiles in $\Sigma^{uc}_j$ as referring to
level $j$ tiles.  The substitution $\phi$ (whereby we think of each
tile of level $j+1$ as the union of four tiles of level $j$) then
defines a map from $\Sigma^{uc}_{j+1}$ to $\Sigma^{uc}_j$ for any
$j\ge 0$, and allows us to define the inverse limit
${\mathop{{\Sigma}}\limits_{\gets}}^{uc}\equiv
\mathop{\lim}\limits_{\longleftarrow} \Sigma^{uc}_j$. 

There is a map from the tiling space $X_{x}$ to
${\mathop{{\Sigma}}\limits_{\gets}}^{uc}$ defined, on $x\in X_{x}$, by
determining precisely how the origin in $x$ sits in each successive
tile of higher and higher level. For many but not all substitution
systems, this map is neither 1-1 nor onto.  A point in the inverse
limit precisely describes the hierarchy of tilies containing the
origin.  In some instances this union of tiles may not cover the
entire plane; if not, it may be that
there is more than one way to extend the tiling to cover the entire
plane (in which case the map is not 1-1), or it may be that there is
no way to extend it (in which case the map is not onto).  For example,
in a chair tiling the tile containing the origin might sit in the
upper left corner of the tile of level 1 of Fig.~2, which sits in the
same position in a corresponding tile of level 2, and so on. The union
of these tiles covers only a quadrant, and this can be extended to a
complete chair tiling of the plane in more than one way.

This problem is absent in substitution systems (such as the Penrose
tilings [AnP]) which ``force the border'', that is, for which
there is some integer $N$ such that, in every tiling, for every tile $\phi^NT$ of level
$N$ the (level 0) tiles which abut this tile are completely determined
by the type $T$.  From this it follows that, for any $M$, the
tiles of level $M$ that abut a given tile of level $N+M$ are also
determined.  An infinite-level tile, therefore, determines a tiling of
the entire plane, and is consistent with a tiling of the entire plane,
even if it does not itself cover the entire plane.  The map from the
space of tilings to the inverse limit is therefore 1-1 and onto.

It is easy to prove, but very useful, that for any substitution system
we can extend the substitution in a simple way to a larger set of
``collared'' tiles producing a tiling space naturally homeomorphic to
the original but now forcing the border. The new set of tiles consists
of multiple marked versions of the original tiles, one for each way
the original tile can be surrounded by tiles in a tiling. (For the
chair one needs 14 versions of each of the original tiles; the
pinwheel requires over 50 versions.)  Any (fixed-orientation)
substitution tiling space can thus be modelled as an inverse limit. 

If the original tiles force the border, we let $\Sigma_x = \Sigma^{uc}$.
If the original tiles do not force the border, let $\Sigma_x$ 
denote the complex constructed from collared tiles. Either way, there
is a natural substitution map from $\Sigma_x$ to itself, and $X_x = 
{\mathop{{\Sigma_x}}\limits_{\gets}}$

The other two kinds of tiling spaces may also be constructed as
inverse limits.  To construct $X_0$ we consider the type {\it but not
orientation} of the tiles of various level that contain the origin.
In the complex $\Sigma_0$ the basic cells are tiles $T$, where each
type appears in only one orientation. As before, if it is possible
for a tile $T_1$ to meet a tile $T_2$ then their common edge is
identified, and it may be necessary to consider collared tiles. As before the
substitution $\phi$ maps $\Sigma_0$ to itself and we consider the
inverse limit.  A point in the inverse limit is a consistent
instruction: the origin sits at such-and-such a point in such-and-such
a tile, which sits in a particular way in a tile of level 1, which
sits in a particular way in a tile of level 2, and so on.  Since we
are using collared tiles (or tiles that force the border), this
prescription defines a unique tiling of the entire plane, up to an
overall rotation, i.e. a point in $X_0$.

The construction of $\Sigma_\phi$ is similar, except that the basic
cells are products $S^1 \times T$, where the $T$ are as in the
construction of $\Sigma_0$.  If a tile $T_1$, rotated by an angle
$\alpha$ (with respect to the standard orientation of $\A$), 
can meet a tile $T_2$, rotated by an angle $\beta$, along a
common edge $e$, then we identify $(\alpha,e) \subset S^1 \times T_1$
with $(\beta,e) \subset S^1 \times T_2$. It then follows that for any
angle $\alpha'$, $(\alpha',e) \subset S^1 \times T_1$ is identified
with $(\alpha'+\beta-\alpha,e) \subset S^1 \times T_2$.  Once again,
substitution maps $\Sigma_\phi$ to itself, and the inverse limit
corresponds to instructions on how to build a tiling around the
origin. Only now the instructions include information on how to orient
each tile, and so defines a tiling uniquely, i.e. a point in $X_\phi$.

If the relative orientation group $G$ is finite then $\Sigma_\phi$,
$\Sigma_x$ and $\Sigma_0$ are all compact simplicial complexes.  This
is clear for $\Sigma_x$ and $\Sigma_0$, as the cells, finite in
number, are simplices that meet along common edges.  The
cells in $\Sigma_\phi$ are not simplices --- they have the topology of
$S^1$ --- but can be divided into $n$ contractible pieces, meeting
along common edges, where $n$ is the order of $G$. 

If the relative orientation group is infinite then $\Sigma_0$ is a
compact simplicial complex but $\Sigma_x$ is not compact (as it
contains an infinite number of cells).  The space $\Sigma_\phi$, while
a compact CW complex, is not necessarily simplicial. 

Finally, these complexes are closely related.  Just as $X_0$ is the
quotient of $X_x$ by $G$ and the quotient of $X_\phi$ by $S^1$,
$\Sigma_0$ is the quotient of $\Sigma_x$ by $G$ and the quotient of
$\Sigma_\phi$ by $S^1$.

\vs.2

\nd {\bf 3. Finite relative orientation group and P-positivity.}

Here we define the ``P-positive cone'' of the top dimensional \Cech
cohomology group for substitution tiling spaces.  In this section we
assume our system has a finite relative orientation group, so the
simplicial complex $\Sigma_x$ is compact.  We also restrict ourselves
to fixed-orientation spaces $X_x$. Later, we will investigate the
role of rotations, and consider order structures on the top
cohomologies of $X_\phi$ and $X_0$ for tiling spaces with arbitrary
relative orientation groups. (We continue to assume our tilings are
2-dimensional, although this construction applies equally well to
other dimensions.)

As discussed above, the space $X_{x}$ is an inverse limit of compact
simplicial complexes $\Sigma_x$ with the substitution (self-)map $\phi$
between them.  One can therefore compute the \Cech cohomology
of $X_{x}$ by taking a direct limit of the \Cech = simplicial cohomology of $\Sigma_{x}$
under the map $\phi^*$.

Recall that the direct limit of an Abelian group under a map $M$ is
the disjoint union of an infinite number of copies of the group,
indexed by the non-negative integers, modulo the equivalence relation
$$ (g,k) \sim (Mg, k+1), \eqno 3)
$$
for all $g$ and all $k$, where $(g,k)$ denotes the element $g$ in the
$k$-th copy of the group.  The direct limit of the cohomology groups
$H^2(\Sigma_{x})$ can be obtained by taking the direct limit of the
cochain groups $C^2(\Sigma_{x})$ then moding out by the image of the
coboundary map $\delta$ from the direct limit of $C^1(\Sigma_{x})$.
(The proof of this is simple diagram chasing.)

Since $\Sigma_{x}$ is a finite simplicial complex, $C_2(\Sigma_{x})$,
the space of simplicial chains, is generated by the tiles themselves
and is isomorphic to $\Z^n$ where $n$ is the number of distinct tiles
in $\Sigma_{x}$. To compute the cohomology, each tile must be given an
orientation. We define this orientation in the following way: Select
an orientation for $\R^2$ and select a tiling $x_0 \in X_x$. The
translations comprise a free $\R^2$-action on $x_0$. Each tile that
appears in $x_0$ thus inherits an orientation from the orientation on
$\R^2$. All such orientations are consistent and are independent of
$x_0$. The projection map from $X_x$ to $\Sigma_x$ defines an
orientation on tiles in $\Sigma_x$. From here on, we refer to this as
the ``positive orientation'' of the tiles in $\Sigma_x$.

If all tiles at all levels are given the positive orientation, the
induced map $\phi_*:C_2(\Sigma_{x}) \to C_2(\Sigma_{x})$ is an $n
\times n$ matrix with only non-negative entries. (The $i,j$ entry in
the matrix is the number of copies of tile type $i$ which appear in a
level 1 tile of type $j$.) By assumption the matrix is primitive, that
is, some power has all of its entries strictly positive.  The cochain
group $C^2(\Sigma_{x})$ is generated by the duals to the tiles.
Relative to this basis, the pullback map $\phi^*$ is described by the
transpose of $\phi_*$. Similar constructions apply to $C^1$, except
that the elements of $C_1$ and $C^1$ are not naturally oriented, so
the substitution matrix on $C^1$ may have negative entries.

Let $\C$ denote the direct limit of the cochains $C^2(\Sigma_{x})$
under the map $\phi^*$, and define $\C_+$ to be the semigroup of
elements $(c,k) \in \C$ which have the property that $(\phi^*)^m c
\in (\Z_+)^n$ for some $m$, where $n$ is the number of distinct tiles
in $\Sigma_{x}$.  It is clear that $\C_+ - \C_+ = \C$.  Since the
matrix $\phi^*$ contains only non-negative entries, it is also clear
that $\C_+ + \C_+ \subset\C_+$, $\C_+ \cap - \C_+ = \{0\}$ and thus
$(\C,\C_+)$ is an ordered group (in fact, a dimension group [Ell]).

Let $H=H^2(X_{x})$ denote the top \Cech cohomology group for $X_{x}$. We
wish to define a positive cone on $H$ by
$$H_+ = \{ [c] \ | \ c \in \C_+\} \eqno 4)$$ 

\nd where elements of $\C$ are equivalent if their difference is
in $\im\, \delta$.
To see that $(H,H_+)$ is an ordered group, we introduce a function
$\mu: \C \to \R$.
\vs.1
\nd {\bf Lemma 1}. {\it There is a homomorphism $\mu : \C \to
\R$ such that $\C_+ -\{0\}= \mu^{-1}(\R_+-\{0\})$.}
\vs.1
\nd {Proof}: 
Since the matrix $\phi^*$ is primitive 
it has a Perron eigenvalue $\lambda>0$ (of multiplicity 1). This is the area
stretching factor of the substitution.  (To see this, note that a
vector whose entries are the areas of the various tiles is an
eigenvector of the substitution matrix.  Since all entries of this
vector are positive this is the Perron eigenvector, and the
corresponding eigenvalue, which is the area stretching factor, must
be the Perron eigenvalue. This eigenvalue is the square of the linear
stretching factor $c(\phi)$.)  Let $r$ denote a left eigenvector of
$\phi^*$.  Define the map $\mu : \C \to \R$ by $(v,k) \mapsto
\lambda^{-k} r v$.

The map is well-defined since 
$$\mu (\phi^* v,k+1) = \lambda^{-k-1} r\phi^* v = 
\lambda^{-k} r v. \eqno 5)$$
Suppose $(v,k) \in \C_+ - \{0\}$ then $(\phi^*)^m v \in (\Z_+)^n$
for some $m$. Since $\phi^*$ is primitive, we may choose $m$ such
that all entries of $(\phi^*)^m v $ are positive. Since all entries
of $r$ are non-negative,
$$\mu(v,k) = \lambda^{-k} rv = \lambda^{-k-m} r (\phi^*)^m v >0. \eqno 6)$$

To prove the reverse direction, take a basis of $\F^n$ which includes
$w_0$, the right eigenvector for $\lambda$, and $(n-1)$ vectors $\{
w_1,w_2,\ldots,w_{n-1}\}$, each $w_i \in \ke (I \lambda_i - \phi^*)^l$ for
some $l$ and $\lambda_i\ne \lambda$. 
We can and do choose $w_0$ with all its entries positive.
Now $r w_i = 0$ since
$$0 = r (I \lambda_i - \phi^*)^l w_i = (\lambda_i - \lambda)^l r w_i. \eqno 7)$$
Therefore if $\mu(v,k) >0$ then writing $v$ in the basis
$\{w_0,w_1,\ldots, w_{n-1}\}$ the coefficient $c_0$ of $w_0$ must be
positive. Since $\lambda$ exceeds the modulus of all other
eigenvalues, as $l \to \infty$, we have $|| \lambda^{-l} (\phi^*)^l v - c_0
w_0|| \to 0$. This implies $(\phi^*)^l v \in (\Z_+)^n$ for some $l$.\qed

Next we show that the map $\mu:\C \to \R$ induces a map from
$H$ to $\R$; since $H = \C/\im \ \delta$, the following will show that
this map is well-defined.
\vs.1
\nd {\bf Lemma 2}. $\im\, \delta \subseteq \ke\, \mu$.
\vs.1
\nd {Proof:} 
Let $e$ be an edge and let $f_e \in C^1(\Sigma_x)$ be the 1-cochain
which assigns a value 1 to the edge $e$ and 0 to all other edges. Then
$\delta(f_e)$ is a 2-cochain which acts on the tiles as follows: let
$T$ be an $m$th level tile. Then $\delta(f_e) (T)$ counts the number
of occurences of $e$ on the boundary $T$ with a $\pm$ signs depending
upon the orientation. In particular, $|\delta(f_e) (T)|$ is less than
or equal to the number of occurences of $e$ on the boundary of the
level $m$ tile $T$.

Now assume $\mu (\delta(f_e)) > 0$. Then there is an $l$ such that for
all for any level 0 tile $T$, $(\phi^*)^l \delta (f_e)(T)>0$.  Furthermore,
for fixed $T$ this quantity scales like $\lambda^l$ for sufficiently
large $l$. This means that for level $l$ tiles there are at least
$O(\lambda^l)$ occurences of $e$ on the boundary. However, $\lambda$
is the area stretching factor, so there can only be
$O(\lambda^{l/2})$ edges on the boundary of a level $l$ tile (since
the perimeter scales like $\lambda^{l/2}$ and there is a minimum
length to each edge).  This is a contradiction.

A similar argument shows that $\mu(\delta(f_e))$ cannot be negative.
Therefore the image of $\delta$ is a subgroup of the kernel of $\mu$
and $\C_+ \cap \im\, \delta = \{0\}$.\qed

Abusing notation, define $\mu:H \to \R$ by $\mu([x]) = \mu(x)$. Letting
$H_+$ be the pre-image under $\mu$ of $\R_+$, $(H,H_+)$ is 
an ordered group since 

\item{(1)} $H_+ + H_+ \subset H_+$

\item{(2)} $H_+ \cap - H_+ = \{0\}$

\item{(3)} $H_+ - H_+ = H$ (this follows from the similar property of $\C_+$).

\vs.2

\nd {\bf 4. Some relations between $X_x$, $X_0$ and $X_\phi$.}

Before we prove the invariance of P-positivity and extend to arbitrary
relative orientation groups, we need to understand the role of
rotations in tilings.  Let $G$ be the relative orientation group of a
tiling and let $R$ be the associated group ring. 

\vs.1

\nd {\it Finite relative orientation groups.} We begin by assuming
$G$ is finite.  As previously indicated, $G$, and therefore $R$, acts
naturally on the simplicial (co)chain groups $C^i(\Sigma_{x})$ and
$C_i(\Sigma_{x})$, giving them the structure of modules over $R$.  A
finite index subgroup of each module can be written as a direct sum of
irreducible representations of $G$.  The substitution maps are
equivariant and so are module homomorphisms, mapping each irreducible
representation to itself. The (co)boundary maps are also equivariant,
and so map a representation in $C^{i}$ (or $C_i$) to the corresponding
representation in $C^{i+1}$ (or $C_{i-1}$).  Thus the integer
(co)homology groups, and in particular the top cohomology, are
themselves finite extensions of finite quotients of direct sums of
irreducible representations.

The need for finite extensions and finite quotients comes about as
follows.  When a compact group acts on $\R^n$, $\R^n$ splits up as the
direct sum of irreducible representations of that group.  However, the
same observation does not apply to $\Z^n$.  For example, if $\Z_2$
acts on $\Z^2$ by permuting the two coordinates, the irreducible
representations are all multiples of $(1,1)$ and all multiples of
$(1,-1)$. The direct sum of these two representations gives all
elements $(a,b)$ with $a+b$ even --- an index 2 subgroup of $\Z^2$.
Similarly, finite index subgroups of the kernels and images of the
(co)boundary maps are direct sums of representations of $G$.
Computing (co)homology one representation at a time is equivalent to
taking a finite index subgroup of the kernel of $\delta$ and modding
out by a finite index subgroup of the image of $\delta$.  Correcting
for our taking too small a numerator means doing a finite extension.
Correcting for the denominator involves taking a finite quotient.

Since the substitution map is equivariant, rotating an eigenvector
gives another eigenvector with the same eigenvalue.  Since the Perron
eigenvectors are unique, with purely positive entries, they must be
rotationally invariant and therefore belong to the trivial
representation. Since the left Perron eigenvector is invariant, $\mu$
of a cochain is also rotationally invariant.  But that implies that
$\mu$ acting on each nontrivial representation is identically
zero. 
  
For example in the uncollared chair complex there is just one tile,
appearing in four different orientations, and $C^2 \sim R=
\Z[t]/(t^4-1)$ where $t$ represents rotation by $\pi/2$.  The
substitution map is given by the $1 \times 1$ matrix $2+t+t^3$. In the
trivial representation we set $t=1$, so our matrix has Perron eigenvalue 4,
with the Perron eigenvector being any multiple of an invariant element of
$C^2$, {\it i.e.} of $(1+t+t^2+t^3)$.

This suggests a streamlined computation of the Perron eigenvalue and
Perron eigenvector.  Instead of working on all of $\Sigma_x$, work on
the simplicial complex $\Sigma_0=\Sigma_{x}/G$, whose (co)chain groups
are precisely the trivial representations of $G$ within the (co)chain
groups of $\Sigma_{x}$.  $C^2(\Sigma_{x}) \otimes_{\Z} \R$ equals
$C^2(\Sigma_0) \otimes_{\Z} \R$ plus the other representations that appear in
$C^2(\Sigma_{x})$.  The Perron eigenvalue on $C^2(\Sigma_{x}) \otimes \R$ is
precisely the Perron eigenvalue on $C^2(\Sigma_0) \otimes \R$, and the (left and
right) Perron eigenvectors on $C^2(\Sigma_{x}) \otimes \R$ come from Perron
eigenvectors on $C^2(\Sigma_0) \otimes \R$. 

Since the (co)chain groups of $\Sigma_0$ are isomorphic to the
invariant parts of the (co)chain groups of $\Sigma_{x}$, and since the
(co)boundary and substitution maps are equivariant, we have
\vs.1
\nd {\bf Theorem 2.} {\it Up to finite extensions,
the (co)homology of $\Sigma_0$ is
equal to the rotationally invariant part of the (co)homology of
$\Sigma_{x}$, and the cohomology of $X_0$ is isomorphic
to the rotationally invariant part of the cohomology of $X_x$.}
\vs.1
We cannot make any statements about the homologies of $X_0$ and $X_x$
because homology, unlike cohomology, does not behave well under
inverse limits.

Now $X_\phi$ is a circle fibration over $X_0$, with singular fibers
corresponding to tilings that are symmetric about the origin.  (These
singular fibers have finite multiplicity, since any given tiling can
only admit discrete rotational symmetry about the origin.)  The
existence of these singular fibers means that computing
$\pi_1(X_\phi)$, $H^1(X_\phi)$ and $H^2(X_\phi)$ from corresponding
data on $X_0$ can be complicated. However, computing the top
cohomology $H^3(X_\phi)$ is quite easy.  From the spectral sequence of
the fibration, it follows immediately that $H^3(X_\phi)$ is isomorphic
to $H^2(X_0)$ [BoT].  From this, and from the fact that $\Sigma_\phi$
is a circle fibration over $\Sigma_0$, we have
\vs.1
\nd {\bf Corollary 1}. {\it Up to finite extensions, the top cohomology of 
$\Sigma_\phi$ is equal to the rotationally invariant part of the top
cohomology of $\Sigma_{x}$. Up to finite extensions, the top cohomology of
$X_\phi$ is isomorphic to the rotationally invariant part of the top
cohomology of $X_x$.}
\vs.1

Note also that the top-dimensional substitution matrix for
$\Sigma_\phi$ is identical to the top-dimensional substitution matrix
for $\Sigma_0$ --- the $S^1$ factor just comes along for the ride. We
may therefore compute our order structure on whichever of the three
spaces $\Sigma_{x}$, $\Sigma_0$ or $\Sigma_\phi$ is most convenient.

To illustrate these principles, we compute the cohomology of the
Penrose tilings.  Here the relative orientation group is $\Z_{10}$,
with group ring $R=\Z[t]/(t^{10}-1)$, where $t$ represents rotation by
$2\pi/10$.  There are 4 irreducible representations of $G$,
corresponding to the factorization $t^{10}-1 =
(t-1)(t+1)(t^4+t^3+t^2+t+1)(t^4-t^3+t^2-t+1)$.  The first two are one
dimensional and $t$ acts by multiplication by $\pm 1$. The other two
representations are 4 dimensional and correspond algebraically to
$\Z[t]/(t^4+t^3+t^2+t+1)$ and $\Z[t]/(t^4-t^3+t^2-t+1)$.

In this tiling, there are 4 types of tiles, whose geometry and
substitutions are shown in Fig.~9.  Similarly, there are four kinds of
edges, each in 10 orientations, and four kinds of vertices, which we
label $\alpha, \beta, \gamma, \delta$, with the relations
$\alpha=t\beta$, $\beta= t\alpha$, $\gamma=t\delta$ and
$\delta=t\gamma$.  The chain groups of $\Sigma_x$ are therefore $C_2 =
C_1 = R^4$ and $C_0 = (\Z[t]/(t^2-1))^2$.  The boundary maps are given
by the matrices
$$ \partial_1= \pmatrix{1-t & -1 & -t & -1 \cr 0 & 1 & 1 & t}; \qquad
\partial_2 = \pmatrix{-1 & t & t^4 & -t^7 \cr -1 & t^9 & -t & t^8 \cr
1 & -t^5 & 0 & 0 \cr 0 & 0 & 1 & -t^5}, \eqno 8)
$$
\nd where we have taken $(\alpha, \gamma)$ as our basis for $C_0$. 
The coboundary maps are given by 
the transposes of these matrices, with $t$ replaced
throughout by $t^{-1}$. The ranks of these matrices are easily computed by row reduction, 
one representation at a time:

\item{(1)} If $t=1$, then $\partial_2^*$ has rank 2 and 
$\partial_1^*$ has rank 1, so $H^0(\Sigma_x)$ and $H^1(\Sigma_x)$
contain one copy of this representation while $H^2(\Sigma_x)$ contains
two.

\item{(2)} If $t=-1$, then $\partial_2^*$ and $\partial_1^*$ have rank 
two, so $H^2(\Sigma_x)$ contains two copies of this representation
while $H^1(\Sigma_x)$ and $H^0(\Sigma_x)$ contain none.

\item{(3)} If $t^4+t^3+t^2+t+1=0$, then $\partial_1^*$ is zero, since 
$C_0$ (also $C^0$) does not contain this representation.
$\partial_2^*$ has rank 4, so all cohomologies are trivial.

\item{(4)} If $t^4-t^3+t^2-t+1=0$, then again $\partial_1^*$ is trivial. 
Now $\partial_2^*$ has rank 3, so $H^1(\Sigma_x)$ and $H^2(\Sigma_x)$
each contain this representation once.

Combining these results, we have that, up to finite extensions,
$H^0(\Sigma_x) = \Z$, $H^1(\Sigma_x) = \Z[t]/(t-1) \oplus
\Z[t]/(t^4-t^3+t^2-t+1)$ and $H^2(\Sigma_x) = (\Z[t]/(t-1))^2 \oplus
(\Z[t]/(t+1))^2 \oplus \Z[t]/(t^4-t^3+t^2-t+1)$.

Next we look at the substitution matrices, which in dimensions 2 and 1 are
$$
\phi_2= \pmatrix{t^7 & 0 & 0 & t^4 \cr 0 & t^3 & t^6 & 0 \cr t^3 & 0 & t^4 & 1 
\cr 0 & t^7 & 1 & t^6}; \qquad \phi_1=\pmatrix{0 & 0 & 0 & t^8 \cr
t^4 & 0 & -t^7 & 0 \cr -t^7 & 0 & 0 & 0 \cr 0 & -t^3 & 0 & -t^5}
\eqno 9)
$$
These matrices (and their pullbacks) are easily seen to be invertible
in all representations, and so induce isomorphisms in cohomology.
Therefore $H^*(X_x)=H^*(\Sigma_x)$.

We can also compute the cohomology of $\Sigma_0$ and $X_0$. The
complex $\Sigma_0$ has four faces, four edges and two vertices.  The
boundary maps can be read off from 8) by replacing $t$ throughout by
1, and the substitution maps are similarly obtained from 9). In other
words, the computation for $\Sigma_0$ and $X_0$ is precisely the same
as the $t=1$ part of the computation for $\Sigma_x$ and $X_x$, with
the result that $H^2(X_0)=H^2(\Sigma_0)=\Z^2$
\vfill 
\eject
\hbox{}
\vbox{\epsfig 1\hsize; 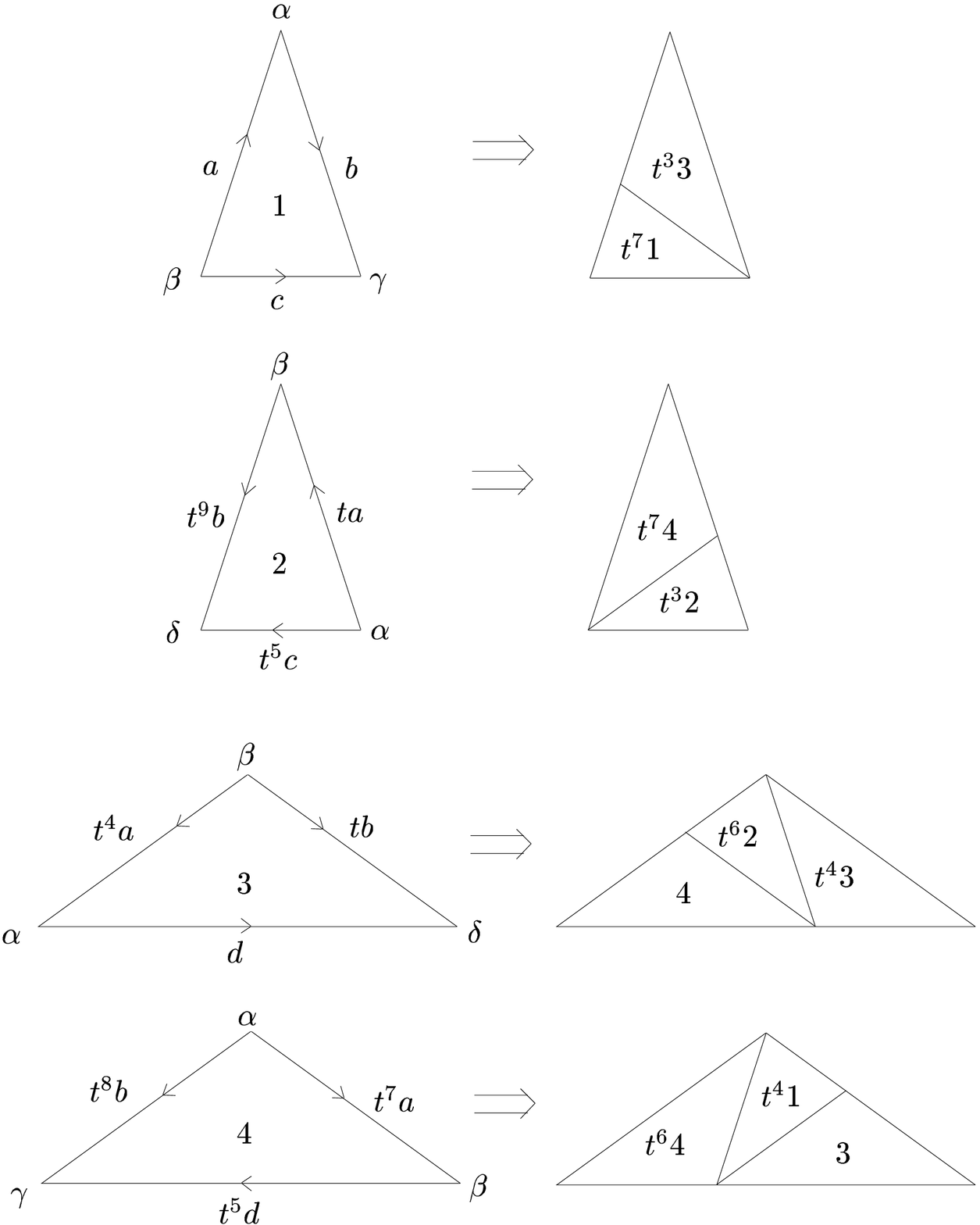;}
\vs.3
\centerline{Figure 9. The Penrose substitution in detail}
\vfill\eject
\nd  and
$H^1(X_0)=H^0(X_0)=H^1(\Sigma_0)=H^0(\Sigma_0)=\Z$.  We therefore also
obtain $H^3(X_\phi)=\Z^2$.

\vs.1
\nd {\it Arbitrary relative orientation groups.}
We have seen that the spaces $\Sigma_x$ and $X_x$ are generally more
complicated than the spaces $\Sigma_0$ and $X_0$, or $\Sigma_\phi$ and
$X_\phi$.  In the case of finite relative orientation group this is a
virtue, as the cohomology of $X_x$ gives us information beyond that of
$X_\phi$. In the case of infinite relative orientation group, however, the
topologies of $\Sigma_x$ and $X_x$ are too rich, and we shall see that
the cohomologies (when nontrivial) are infinitely generated.  This
makes them too big to handle with the formalism of primitive matrices
and Perron eigenvectors.  Instead, we are forced to work on $\Sigma_0$
and $X_0$ (or $\Sigma_\phi$ and $X_\phi$) where the cohomology groups
are finitely generated and we can define an order structure.  This
structure then pulls back to an ordering of the rotationally-invariant
elements of the top cohomology of $X_x$.
 
If a tiling space has an infinite relative orientation group,
then the simplicial complex $\Sigma_x$ is noncompact, but $\Sigma_0$
is still compact, and has a finite substitution matrix.
In that case, the top homology of $\Sigma_x$, if nontrivial, is infinitely
generated.  This is seen as follows.  Let $C$ be a boundaryless
2-chain on $\Sigma_x$.  Since $C$ is a finite collection of tiles,
only a finite number of orientations appear in $C$.  By applying an
appropriate rotation $g$, we obtain another boundaryless 2-chain
$C'=g(C)$ whose tiles appear in a completely different set of
orientations from those of $C$.  Similarly, $C''=g(C')$ is different
from $C'$. If the element $g$ has infinite order, we can generate an
infinite collection of boundaryless chains, all linearly independent
in $C_2(\Sigma_x)$.  Thus $H_2(\Sigma_x)$ is infinitely generated.  A
similar argument, involving the dual to a tile in $C$, shows that
$H^2$ is also infinitely generated.  Furthermore, since the dual to a
direct sum is a direct product, $H^2$ contains elements represented by
infinite sums of duals of tiles.

Finding the left Perron eigenvector is also a problem.  The left
Perron eigenvector of $\phi^*$ is the same as a right Perron
eigenvector of $\phi$.  This must be rotationally invariant, and is
therefore a sum of an infinite number of tiles.  Pairing this infinite
sum with an element of $C^2$, which itself has an infinite number of
terms, gives an infinite sum of real numbers that need not converge.
Thus the $\mu$ map on $H^2(\Sigma_x)$, or on $C^2(\Sigma_x)$, is not
well defined. 

Working on $\Sigma_0$, however, we have no problem. The space $\Sigma_0$ is a
compact simplicial complex, with finitely generated (co)chain
complexes, and finitely generated (co)homology.  The substitution
matrix $\phi^*$ is finite and primitive, and has a Perron eigenvalue
and Perron eigenvector exactly as before.  On $H^2(\Sigma_0)$, and on
the limit $H^2(X_0)$ ($=H^3(X_\phi)$) the order structure defined
above works exactly as in the case of finite relative orientation
group.
\vs.2

\noindent {\bf 5. The invariance of positivity.}

Having defined the P-positive cone in the top \Cech cohomology group
of our tiling spaces, we wish to show that this positive cone is
invariant (up to an overall sign) under homeomorphisms. (This is not
immediate, since the inverse limit structure that we used to define
the positive cone is {\it not} invariant.)

In general, the \Cech cohomology of a space is the direct limit of
\Cech cohomologies of open covers where the open covers generate 
the topology of the space (see [BoT] for details).  In our situation
we begin by considering a triangulation of $\Sigma_n$ and the open
cover of $\Sigma_n$ given by ``open stars'' of the vertices of the
triangulation. By the open star ($\star(v)$) of a vertex $v$ we mean
the collection of all simplices that touch $v$, but not including the
opposite boundary simplex. We associate to this open cover a
simplicial complex called the ``nerve'' of the open cover. The nerve
of an open cover is constructed by including one vertex for each open
set in the cover, then including a $k$-simplex between any $k$-vertices
such that the intersection of the corresponding $k$ sets in the open
cover is nonempty. The \Cech cohomology of the cover is defined
as the simplicial cohomology of the nerve of that cover.  Since every
vertex of the triangulation lies in exactly one open star and every
$k$-simplex lies in exactly $k$ open stars, the nerve of the open star
cover of $\Sigma_n$ is identical to the triangulation of $\Sigma_n$ as
a simplicial complex. Thus the \Cech cohomology of the open star cover
of $\Sigma_n$ is {\it canonically} identified with the simplicial
cohomology of the triangulation of $\Sigma_n$. (This argument shows
that simplicial and \Cech cohomology agree for all simplicial
complexes.)

Consider triangulations of $\Sigma_n$, generated recursively, so that
each triangle in the triangulation of $\Sigma_n$ maps by $\phi$ to a
single triangle in the first barycentric subdivision of the
triangulation of $\Sigma_{n-1}$. This triangulation of $\Sigma_n$
corresponds to an open cover of $\Sigma_n$, which pulls back by the
natural projection $\pi_n: X_x \to \Sigma_n$ to give an open cover
${\cal A}_n$ of $X_x$.  That is, an open set in ${\cal A}_n$ is the
set of tilings in which the origin sits in a specified neighborhood in
a tile of level $n$.  These neighborhoods, as defined by the
triangulations of $\Sigma_n$, are increasingly fine as $n \to \infty$
(even after accounting for the expansion due to $\phi$), so the
sequence of open covers ${\cal A}_n$ generates the topology of $X_x$.

The \Cech cohomology of $X_x$ is the direct limit of the \Cech
cohomologies of the covers ${\cal A}_n$ (hence the direct limit of the
simplicial cohomologies of $\Sigma_n$) with a bonding map defined as
follows. From their construction, ${\cal A}_{n+1}$ refines ${\cal
A}_n$. Each element of ${\cal A}_n$ (resp. ${\cal A}_{n+1}$) is a set
of the form $\pi_n^{-1}(\star(v))$ where $v$ is a vertex in the
triangulation of $\Sigma_n$ (resp. $\Sigma_{n+1}$). We define a vertex
map from $\Sigma_{n+1}$ to $\Sigma_n$ by mapping each vertex $u$ in
the triangulation of $\Sigma_{n+1}$ to a vertex $v$ in the
triangulation of $\Sigma_n$ where $\pi_{n+1}^{-1}(\star(u)) \subset
\pi_n^{-1}(\star(v))$. After extending this vertex map linearly to 
the interior of higher level simplices, the result is our bonding map
$f_n: \Sigma_{n+1} \to \Sigma_n$. Although the vertex map, and
therefore $f_n$, may not be uniquely defined, the induced chain map
$f_n^*$ on cohomology is always the same [BoT].

To define positivity recall that we used an orientation in $\R^2$ to
define an orientation of all 2-simplices in $\Sigma_n$.  Assume
$\alpha$ and $\beta$ are two such oriented 2-simplicies which both
occur within the same tile $T \in \Sigma_n$. Then the duals to
$\alpha$ and $\beta$ are cohomologous. Any positive sum of these
cohomology classes was declared to be P-positive in section 3.

To show that this notion of P-positivity is invariant we must
consider the isomorphism of \Cech cohomology induced by a
homeomorphism $h$. Let $\{\Sigma_n\}$ and $\{\Sigma_n'\}$ be sequences
of complexes corresponding to the tiling spaces and let $\{{\cal
A}_n\}$ and $\{{\cal A}_n'\}$ be sequences of open covers, of $X_x$
and $X_x'$, resp., corresponding to triangulations of the
complexes. Let ${\cal B}_n = h^{-1}({\cal A}_n')$. From their
construction, for every $n$ there is an $m$ such that ${\cal A}_m$
refines ${\cal B}_n$. Thus we may define the bonding map $f_{{\cal B}_n
{\cal A}_m}: \Sigma_m \to \Sigma_n'$ by mapping each vertex $u$ to a
vertex $v$ where $(\pi_m)^{-1}(\star(u)) \subset
h^{-1}({\pi_n'}^{-1}(\star(v)))$.  By taking limits we obtain an isomorphism
between the \Cech cohomology of $X_x$ and the \Cech cohomology of
$X'_x$.

\vs.1
\nd {\bf Theorem 3.} {\it If $X_x$ and $X_x'$ are two fixed-orientation 
tiling spaces with finite relative orientation groups, and $h:X_x \to
X_x'$ is a homeomorphism, then the isomorphism $h^* :H^2(X_x') \to
H^2(X_x)$ preserves the P-positive cones, up to an overall sign.}

\vs.1
\nd Proof:
To see that the P-positive cone is preserved by $h^*$, first
consider whether the map $h:X_x \to X_x'$ is orientation preserving or
reversing on each path-component. 
That is, fix a tiling $x_0 \in X_x$ and $h(x_0) \in
X_x'$. There is a 1:1 map $j:\R^2 \to \R^2$ defined by $h(g(x_0)) =
j(g) h(x_0)$ where $g \in \R^2$ acts by translation. We say
$h$ is orientation preserving or reversing based on whether $j$ is
orientation preserving or reversing as a homeomorphism of
$\R^2$. Since the orbit of each point is dense, and since the action 
of $\R^2$ is continuous, whether $h$ is
orientation preserving or not is independent of the choice of
$x_0$. Without loss of generality, assume that $h$ is orientation
preserving.

The goal is to show 
that every P-positive cohomology class 
on $X'$ pulls back to a P-positive cohomology class on $X$.
But the P-positive classes on $X'$ are precisely the classes represented by
positive linear combinations
of basis cochains on $\Sigma_n'$, with $n$ arbitrarily large.  So we
take a basic elementary cochain on $\Sigma_n'$ (with $n$ large), 
pull it back to $X$, and try 
to show that it's positive there.  Without loss of generality we can assume 
that this cochain, which we call $1_\Delta$, is dual 
to a triangle $\Delta$ sitting in the middle
of one of the tiles of $\Sigma_n'$.  We then find an $m$ large enough
so that the open cover induced by the partition of $\Sigma_m$ is a refinement
of the open cover induced by the partition of $\Sigma_n'$.  This refinement
gives a vertex map $\rho$ from the nerve of the cover of $\Sigma_m$ (that is, 
$\Sigma_m$ itself) to the nerve of the cover of $\Sigma_n'$ (that is, 
$\Sigma_n'$ itself).  We must show that $\alpha= \rho^* 1_\Delta$ represents a 
P-positive cohomology class.

Now $\alpha$ itself is not necessarily a positive linear combination
of elementary cochains.  There will be some triangles $T$ on $\Sigma_m$
where $\alpha$ evaluates to $+1$ or $-1$, and others where it evaluates
to 0.  We call the triangles where it evaluates to $\pm 1$ {\it essential}.  
Our task is to show that the positive essential triangles  outnumber
the negative ones in the following sense: there exists an $N>m$ such that
every complete tile in $\Sigma_N$, mapped to $\Sigma_m$ and paired with
$\alpha$, yields a positive number.  Put another way, $\alpha$ evaluates
positively on the image (in $\Sigma_m$) of every tile of 
level $N$ in a particular tiling in $X$. 

Here we are using the identification between a particular tiling in $X$
(viewed as $\R^2$
with appropriate markers) and a path-component in $X$. In this 
identification, a point on a tiling corresponds to a translate of that tiling 
with that particular point situated at the origin.  Call this space $P$ (for
path-connected) and let $P'$ be the corresponding path-component of $X'$. 
The triangulation of $\Sigma_m$ induces a triangulation of $P$, with each
vertex, edge, or triangle of $P$ mapping by $\pi_m$ to a corresponding 
vertex, edge, or triangle of $\Sigma_m$.
There are two metrics on
$P$, namely the Euclidean metric and the metric inherited from $X$.
For short-to-moderate distances, these metrics behave essentially identically:
\vs.1
\noindent {\bf Lemma 3}. {\it There exist constants $\epsilon_1$, $\epsilon_2$ 
and $M$, with $0<\epsilon_1 <M$ and $0< \epsilon_2$, 
such that any two points in the same path component with 
Euclidean distance in $(\epsilon_1, M)$ have tiling-space distance greater
than $\epsilon_2$.}  
\vs.1

\noindent Sketch of Proof: For our fixed set of polygonal tiles let $K$ be a
lower bound on the distance between the midpoint of any edge of any
tile in any tiling and any other edge in that tiling, and let $\alpha$
be a lower bound of the interior angles of any of the tiles. Let $N$
be such that any open ball of radius $N$ in any tiling contains at
least one tile and all the tiles touching it. The result then follows
if $M < \min \{K/2,\ \sin(\alpha/2)\}$ and $\epsilon_2 < \epsilon_1/N$.\qed

\vs.1
\noindent {\bf Corollary 2}. {\it A path-connected set of tiling-space diameter
$\epsilon_2$ or less has Euclidean diameter $\epsilon_1$ or less.} 

We return to the proof of Theorem 3. 

Pick $\epsilon_{1,2}$ for $X$ as in the lemma.  Since our
homeomorphism $h$ is between compact spaces, $h^{-1}$ is uniformly continuous,
so there is a $\delta$ for which all sets of (tiling-space) diamater
$\delta$ or less in $X'$ get mapped to sets of (tiling-space) diameter
$\epsilon_2$ or less in $X$. Then pick $n$ large enough that we can choose
our triangle $\Delta$ to have Euclidean 
diameter $\delta$ or less.  This implies that each path-connected component
(PCC) of $h^{-1}(\Delta)$ in $X$ will have Euclidean diameter 
$\epsilon_1$ or less.

Our proof of positivity now proceeds in four steps:

1. Show that preimages in $P$ of essential triangles appear in
   disjoint and localized clusters, with each cluster associated to a
   specific PCC of $h^{-1}(\Delta)$.  
   
2. Show that $\alpha$ applied to each complete cluster gives exactly $+1$. 

3. Show that the number of complete clusters in a tile (of any level) is bounded
   from below by a multiple of the area of the tile. 

4. Show that $\alpha$, applied to triangles whose clusters are only
   partially in the tile, is bounded in magnitude from above by a
   multiple of the perimeter of the tile.  (This contribution
   need not be positive.)

Since volume scales faster than perimeter, $\alpha$ then evaluates 
positively on all tiles of sufficiently high level.

Step 1. If $T$ is an essential triangle, then its 3 vertices get sent
by the bonding map to the 3 vertices of $\Delta$.  But this means that
the open star of each of the vertices of $T$ is contained in $h^{-1}$
of the open star of the corresponding vertex of $\Delta$. But then the
intersection of these three open stars (namely $T$ itself) is
contained in $h^{-1}$ of the intersection of the open stars of the
vertices of $\Delta$ (namely $\Delta$ itself).  Thus the open set $T$ is
contained in the open set $h^{-1}(\Delta)$, and every PCC of $T$
must lie in a PCC of $h^{-1}(\Delta)$.  

We define a cluster to be all the PCCs of essential triangles within a
specific PCC of $h^{-1}(\Delta)$.  The PCCs of $\Delta$ in $P'$ are
disjoint (since $\Delta$ was chosen to lie in the interior of a tile
in $\Sigma_n'$) and have Euclidean diameter $\delta$ or less. This
implies that the PCCs of $h^{-1}(\Delta)$ in $P$, and therefore the
clusters, are disjoint and have Euclidean diameter $\epsilon_1$ or
less.

Step 2.  Let $D$ be a PCC of $h^{-1}(\Delta)$, and let $\Gamma$ 
be the sum of all the triangles that intersect $D$. 
$\Gamma$ is one of our clusters, 
plus a number of components of inessential triangles. 
We must show that $\alpha(\Gamma)=+1$.  We do by showing that, as 
a simplicial chain, $\rho(\pi_m(\Gamma))=\Delta$. 

Let $\{v_i\}$ be the collection of all vertices in $\Gamma$. This is
precisely the set of vertices whose stars contain points of $D$. 
If $v_i$ is one of the vertices in this collection then $\rho(\pi_m(v_i))$
must be one of the vertices, $d_1, d_2, d_3$,
of $\Delta$. This follows since
$(\star(v_i))$ contains a point in $D$, and the image of this point 
lies in exactly three open sets in the open cover of $\Sigma_n'$, namely
the stars of the three vertices of $\Delta$.  Since the cover induced
by the triangulation of $\Sigma_m$ is a refinement of the cover induced by
the triangulation of $\Sigma_n'$, the open set corresponding to $\pi_m(v_i)$
is contained in one (or more) or these three open sets, and in no others.
Thus if $\gamma$ is a 2-simplex in
$\Gamma$ then $\rho(\pi_n(\gamma))$ is either $\pm
\Delta$ or zero. Therefore, $\rho (\pi_m(\Gamma)) = k \Delta$  for some
integer $k$. We will show that $k=1$ by showing that the boundary of
$\Gamma$ maps to $+1$ times the boundary of $\Delta$.

The boundary of  $\Gamma$  is a closed loop in $X_x$. 
Let $x_1, x_2, x_3$ be the three vertices of $D$ (so that 
$\{h(x_i)\}$ are the 3 vertices of a PCC of $\Delta$). We may
choose 3 distinguished vertices $v_1, v_2, v_3$ on the boundary of
$\Gamma$ such that $\star(v_i)$ contains
$d_i$ for $i=1,2,3$. Because $x_i$ is contained in a unique open set
in $\B_n$ 
(namely that corresponding to the $d_i$),
we must have $\rho(\pi_n(v_i)) = d_i$ for all $i$.  
Now let $p_{12}$ denote the path from
$x_1$ to $x_2$ along the boundary of $D$. 
There is a portion of the boundary
of $\widetilde \Gamma$ which forms a path from $v_1$ to $v_2$ with the
property that the star of each vertex along this path 
intersects $p_{12}$. Since the only elements
of ${\cal B}_m$ that can contain points in $p_{12}$ are
the open sets corresponding to $d_1$ and $d_2$, the simplicial map
$\rho \circ \pi_m$ applied to this portion of the boundary of 
$\Gamma$ gives a path that begins at $d_1$, ends at $d_2$ and is
completely contained in the edge between $d_1$ and $d_2$.  The net
output of the chain map $\rho \circ \pi_m$ applied to this path is
therefore $+1$ times the edge from $d_1$ to $d_2$. Applying the same
reasoning around the entire boundary of $\Gamma$ and
recalling that the map $h$ is orientation-preserving, we see that
$\rho \circ \pi_m$ applied to boundary of
$\Gamma$ is exactly $\partial \Delta$, as required.

Step 3. Since each preimage of an essential triangle has fixed
Euclidean area, and since each PCC of $h^{-1}(\Delta)$ has diameter
$\epsilon_1$ or less, there is an upper bound to how many essential
triangles can fit in a cluster.  The number of clusters is then
bounded from below by the number of PCCs of essential triangles,
divided by the maximum number of triangles per cluster. But each
essential triangle appears with prescribed positive density, so the
number of clusters must scale with the area.  Of these we must
subtract off the partial clusters near the boundary, but this scales
like the perimeter.  See Step 4, below.

Step 4. Of the previously considered clusters, we must give special
treatment to the ones within a distance $\epsilon_1$ of the boundary
of the tile.  These may lie only partially within the tile,
and the contribution of the essential triangles within the tile
may not be positive.  However, the area of this boundary region is at
most $\epsilon_1$ times the perimeter of the tile. Since the area
of each triangle is bounded from below, the number of essential
triangles in this region is bounded by a multiple of the perimeter. \qed

\vs.2

\nd {\bf 6. Extension to $X_\phi$ for arbitrary relative orientation groups.}

The previous development was for the cohomology of $X_x$ in the case
where the relative orientation group is finite. We now consider the
\Cech cohomology of $X_\phi$ with no restrictions on the relative
orientation group. 

In the case of arbitrary relative orientation group, P-positivity was
defined for $H^2(X_0)$, using Perron eigenvectors and the substitution matrix, 
at the end of section 4.  The natural isomorphisms between
$H^3(X_\phi)$ and $H^2(X_0)$, and between $H^3(\Sigma_\phi)$ and
$H^2(\Sigma_0)$, then give a consistent definition of P-positivity
for $H^3(X_\phi)$.  This notion of positivity can be rephrased 
directly in the language of \v Cech cohomology.
 
When working with $X_\phi$, our complexes $\Sigma_n$ are comprised of
$S^1 \times$ tiles. Triangulations of $\Sigma_n$ then contain
vertices, edges, triangles and tetrahedra. Correspondingly, our open
covers, by stars of vertices, have single, double, triple and
quadruple intersections. We orient the tetrahedra according
to some global orientation on the 2-dimensional Euclidean group. The
dual to any two 3-simplices within any $S^1 \times$ tile are
cohomologous. We declare that the positive sum of any of the duals to
any of these 3-simplices is positive.  A cohomology class on $X_\phi$
is positive if it can be represented as the pullback from some $\Sigma_n$ of such a 
positive sum.  

\vs.1 
\nd {\bf Theorem 4}. {\it If $X_\phi$ and $X_\phi'$ are two all-orientation 
tiling spaces, and $h:X_\phi \to
X_\phi'$ is a homeomorphism, then the isomorphism $h^* :H^3(X_x') \to
H^3(X_x)$ preserves the P-positive cones, up to an overall sign.}

\noindent Sketch of Proof: 
The proof of Theorem 3
carries over with only minor modification, as sketched below.

We first consider whether $h$ preserves or
reverses the orientations of path components of the tilings spaces --
that is, whether it preserves or reverses the orientations of the
action of the Euclidean group.  Without loss of generality, assume
that $h$ is orientation-preserving.

If $\Delta$ is a 3-simplex arising from the triangulation
of $\Sigma'_n$, with $n$ large, we choose $m$ such that the open cover
induced by the triangulation of $\Sigma_m$ refines the open cover 
induced by the triangulation of $\Sigma'_n$, pulled back to $X_\phi$ by $h^{-1}$. 
We then define the cochain $\alpha=\rho^* 1_\Delta$ on $\Sigma_n$, where
$\rho: \Sigma_m \to \Sigma'_n$ is a vertex corresponding to the refinement.
We call a tetrahedron on $\Sigma$ {\it essential} if $\alpha$ applied to that
tetrahedron is nonzero.  We must show that every tile of sufficiently high
level, crossed with $S^1$, contains more preimages of positive essential
tetrahedra than of negative essential tetrahedra. 

As before we work on a particular path-component $P$ of $X_\phi$. $P$
is naturally identified, not with a single tiling, but with the unit
tangent bundle of a single tiling as follows. A unit vector at a particular
point in the tiling means that the tiling should be translated to put
that point at the origin, and rotated so that the vector points in the
positive $x$-direction. This unit tangent bundle has a natural
Riemannian metric that plays the same role as that played by the
Euclidean metric in the proof of Theorem 3. In particular, analogs
to Lemma 3 and Corollary 2 are easily proven.

The proof then proceeds in 4 steps, exactly as before.  We show that preimages
of essential tetrahedra appear in clusters, that $\alpha$ applied to
$\pi_m$ of each complete cluster is $+1$, that the number of complete
clusters in a tile of high level (crossed with $S^1$) 
scales as the area of the tile, and
that the effect of the incomplete clusters is bounded by a multiple of
the perimeter of the tile.

The only subtlety is in step 2, where previously we had made explicit use
of the geometry of triangles. We define $D$ and $\Gamma$ as
before, and show that each vertex in $\Gamma$ must be mapped by $\rho
\circ \pi_m$ to a vertex of $\Delta$.  Those vertices in $\Gamma$
whose open stars contain a vertex $x_i$ of $D$, must get mapped to the vertex
$d_i$ of $\Delta$.  Those whose open stars intersect the edge of $D$ that runs
from $x_i$ to $x_j$ must get mapped to $d_i$ or $d_j$.  Those that
intersect the face whose corners are $x_i$, $x_j$ and $x_k$ must be
mapped to $d_i$, $d_j$ or $d_k$. These considerations, together with
the fact that $\Gamma$ is contractible, imply that $\partial \Gamma$,
viewed as a chain, is mapped to $+1$ times $\partial \Delta$, and hence that
$\Gamma$ is mapped to exactly $+1$ times $\Delta$.  \qed

\vs.2
\nd {\bf 7. Algebraic consequences.}

Below we examine some of the consequences of isomorphism of ordered
groups obtained as direct limits.  By Theorems 3 and 4 this will show
that certain properties of the area stretching factor are
homeomorphism invariants of $X_x$ for the finite orientation case, and
are homeomorphism invariants of $X_\phi$ in all cases.

Our setting is as follows.  We assume we have Abelian groups $H$ and
$H'$, obtained as direct limits of finitely generated free Abelian
groups by primitive maps $\phi^*$ with Perron eigenvalues $\lambda$
and $\lambda'$.  An element of $H$ will be denoted $[(v,k)]$, meaning
the class of the element $v$ in the $k$-th approximant to $H$.  Of
course, the application we have in mind is for $H$ to be the top \v
Cech cohomology of a tiling space, for the free Abelian groups to be
the free part of the top cohomology of $\Sigma_x$ or $\Sigma_\phi$,
and for the primitive map to be induced by substitution.

\vs.1
\nd {\bf Lemma 4}. {\it We can choose the left Perron eigenvector  
$r = (r_1, r_2, \ldots, r_n)$  of the primitive matrix $\phi^*$
so that each $r_i$ is a polynomial in $\lambda$ with integer coefficients.}
\vs.1
\nd {Proof:} To find $r$, we wish to solve the equation
$$r (I \lambda - \phi^*) = 0. \eqno 11)$$ 
In order to solve this, we perform
column operations until the matrix $\phi^*$ is in
lower triangular form. Throughout this process, each entry of the
matrix will be a quotient of integer polynomials in
$\lambda$. Therefore there is a solution to the equation where each
entry of $r$ is a quotient of integer polynomials in $\lambda$. Thus
by rescaling we may obtain a row eigenvector such that each entry is
an integer polynomial in $\lambda$. \qed

\noindent 
{\it Integer Case:}
In the case where $\lambda$ is an integer, we will take the $r$ from
the previous lemma and scale by a factor of $1/g$ where $g$ is the
greatest common divisor of the entries of $r$.

\noindent 
{\it Noninteger Case:}
In the case where $\lambda$ is not an integer (and therefore
irrational), take the $r$ from the previous lemma and scale by a
factor of $1/\lambda^p$ where $p$ is the maximal power of $\lambda$
that occurs in the entries of $r$. Now the row eigenvector is in 
$\Z[1/\lambda]$. 

We define the map $\mu:H \to \R$ by $\mu[(v,k)] \mapsto \lambda^{-k} r v$ where $r$ is
fixed as above.  In either the integer or noninteger case it is clear that the image of $\mu$ is
a subset of $\Z[1/\lambda]$. Our positive cone $H_+$ is then the preimage
of the positive real numbers.

Now we suppose that we have two such groups, $H$ and $H'$, with maps 
$\mu$ and $\mu'$ from each to $\R$, constructed as above,
with Perron eigenvalues $\lambda$ and $\lambda'$, and suppose that $h:(H,H_+)
\to (H',H_+')$ is an order isomorphism. (We do not assume that $h$ 
respects the direct limit structure, only that it respects the order
structure.)
\vs.1

\nd {\bf Lemma 5}. {\it Let $x,y \in H$ with $\mu(x)$ and $\mu(y)$ both 
positive.  Then} 
$${\mu(x) \over \mu(y)} = {\mu'h(x) \over
\mu' h (y)}. \eqno 12)$$
\vs.1
\nd {Proof:} 
Let $a,b$ be positive integers. Then 
$$\eqalign{
a x - b y \in H_+ 
& \iff \mu(ax-by) \geq 0 \cr
& \iff a \mu(x)-b \mu(y) \geq 0 \cr
& \iff {\mu(x) \over \mu(y)} \geq {b \over a}.} \eqno 13)$$

\nd Since $h$ is an order isomorphism, 
$$\eqalign{ a x - b y \in H_+ 
& \iff a h(x) - b h(y) \in H_+'  \cr
& \iff {\mu' h(x) \over \mu' h(y)} \geq {b \over a}.} \eqno 14)$$

\nd Therefore $b/a \leq \mu(x)/\mu(y)$ if and only if 
$b/a \leq \mu' h(x)/ \mu' h(y)$. Since this
is true for any positive integers $a$ and $b$ it must be the
case that $${\mu'h(x) \over \mu' h(y)} =  {\mu(x) \over\mu(y)}.\qed \eqno 15)$$

\vs.1
\nd {\bf Lemma 6}. 
{\it $\Q [\lambda]$ and $\Q [\lambda']$ are identical as subsets of $\R$.}
\vs.1
\nd {Proof:}
First we claim that there exist $p(x), q(x) \in \Z[x]$ such that
$$\lambda = {p(\lambda') \over q(\lambda')}. \eqno 16)$$

\nd Choose $x = [(v,k)] \in H$ such that $\mu(x) >0$. Let $y = [(v,k+1)]
\in H$.  Then 
$${\mu' h(x) \over \mu' h(y)}= {\mu(x) \over \mu(y)} = \lambda. \eqno
17)$$ 
Since the image of $\mu'$ is a subset of $\Z [1/\lambda']$, the
claim follows.

Now let $f(x)$ be the minimal monic polynomial in $\Z[x]$ for the
algebraic integer $\lambda'$. Since $q(\lambda') \neq 0$, the
polynomial $f(x)$ does not divide $q(x)$. Since $f(x)$ is irreducible,
the number 1 is a greatest common divisor of $f(x)$ and $q(x)$ in the
ring $\Q[x]$.

Since $\Q[x]$ is a principal ideal domain, there exist polynomials
$a(x), b(x) \in \Q[x]$ such that 
$$a(x) q(x) + b(x) f(x) = 1. \eqno 18)$$
This implies
$$a(\lambda') q(\lambda') = 1. \eqno 19)$$ 

\nd Therefore, 
$$\lambda = a(\lambda') q(\lambda') \in \Q[\lambda'].\eqno 20) $$
which implies that as sets of real numbers, 
$$\Q[\lambda] \subseteq \Q[\lambda']. \eqno 21)$$
The other inclusion follows similarly. 
\qed

(Of course, the previous lemma is trivial in the case where
$\lambda$, $\lambda'$ are integers.)
\vs.1
\nd {\bf Lemma 7}. {\it If $\lambda$ is an integer then $\lambda'$ is also 
an integer. Furthermore, $\lambda$ and $\lambda'$ have the same prime
factors.}
\vs.1
\nd {Proof:}
Assume $\lambda$ is an integer and recall $\mu [(v,k)]$ is defined as
$\lambda^{-k} rv$ where $r$ is an integer vector with relatively prime
entries. Therefore there are integers $v_1, v_2, \ldots , v_n$ such that 
$$\sum_{i=1}^n v_i r_i =1. \eqno 22)$$
Define $v$ to be the column vector $\lambda (v_1, v_2, \ldots , v_n)^{\rm t}$. 
Then $\mu[(v,1)] = 1$ and the number 1 is in the image of $\mu$.  

 From the previous lemma it is clear that $\lambda'$ must also be an
integer (it is a Perron eigenvalue which is either an integer or
irrational). Thus for the same reason 1 is in the image of $\mu'$. 

Suppose $y \in H'$ with $\mu'(y) =1$. Let $[(u,k)] = h^{-1}(y)$.  Then
$$\eqalign{
{\mu[(u,k+1)] \over \mu [(u,k)]} & = {\mu' h[(u,k+1)] \over \mu' (y)}\cr
{1 \over \lambda} & = \mu' h[(u,k+1)].} \eqno 23) $$

\nd This implies $1/\lambda \in \Z[1/\lambda']$. Similarly, 
$1/\lambda' \in \Z[1/\lambda]$. Therefore, as subsets of $\R$, 
$$\Z[1/\lambda] = \Z[1/\lambda']. \eqno 24)$$

\nd This is the case if and only if $\lambda$ and $\lambda'$ have the same
prime factors.\qed
\vs.1

Combining these results with Theorems 3 and 4 then yields our main
result, Theorem 1.
\vs.2 \nd
{\bf 8. Conclusion.}

We add here some remarks of a general nature. The basic objects of the
above study are hierarchical tilings of the Euclidean plane, which are
mainly of interest for their unusual symmetry (geometric)
properties. Our analysis is an outgrowth of a traditional dynamical
analysis, in which a space of tilings is equipped with a topological
and Borel structure together with the action of Euclidean motions, the
latter maintaining contact with geometry. It is notable that we have
been ignoring the dynamical action, obtaining topological invariants
of tiling spaces. This is a recent development, but not new to this
paper, appearing as it has in noncommutative approaches. However a
distinctive feature of this work is that our invariant is only
sensible in the top (\v Cech) cohomology, which suggests to us that it
might not be essentially topological; more specifically, it might be
invariant under maps more general than homeomorphisms.

\vs.2 \nd
{\bf Acknowledgements.}\ We gratefully acknowledge useful
discussions with Mike Boyle, Gary Hamrick, Ian Putnam, Mike Starbird and Bob
Williams.
\vfill\eject
\nd {\bf References}
\vs.3 \nd
[AnP]\ J.\ Anderson and I.\ Putnam, Topological invariants for substitution
tilings and their associated C$^\star$-algebras, {\it Erg.\ Thy.\ Dyn.\ Syst.}\ {\bf 18} 
(1998), 509-537.
\vs.1 \nd
[BJV]\ M.\ Barge and B.\ Diamond, A complete invariant for the topology of
substitution tiling spaces, preprint.
\vs.1 \nd
[BoT]\ R.\ Bott and L.\ Tu, {\it Differential forms in algebraic topology},
Graduate Texts in Mathematics, {\bf 82}, Springer-Verlag, New York, 1982.
\vs.1 \nd
[Con]\  A.\ Connes, {\it Noncommutative Geometry}, Academic Press, San Diego, 1994.
\vs.1 \nd
[Ell]\ G.\ Elliott, On the classification of inductive limits of
sequences of semisimple finite dimensional algebras, {\it J.\ Algebra} {\bf 38} 
(1976), 29-44.
\vs.1 \nd
[Goo]\ C.\ Goodman-Strauss, Matching rules and substitution tilings,
{\it Annals of Math.}, {\bf 147} (1998), 181-223.
\vs.1 \nd 
[GPS]\ T.\ Giordano, I.\ Putnam and C.\ Skau, Topological orbit equivalence and 
C$^*$-crossed products, {\it J.\ reine angew.\ Math.} {\bf 469} (1995), 51-111.
\vs.1 \nd
[HPS]\ R.\ Herman, I.\ Putnam and C.\ Skau, Ordered Bratteli diagrams, 
dimension groups, and topological dynamics, {\it Intern.\ J.\ Math.} {\bf 3} (1992),
827-864.
\vs.1 \nd
[Ke1]\ J.\ Kellendonk, Non-commutative geometry of tilings and gap 
labelling, {\it Rev.\  Math.\ Phys.}, {\bf 7} (1995), 1133-1180.
\vs.1 \nd
[Ke2] J.\ Kellendonk, Topological equivalence of tilings, {\it
J.\ Math. Phys.} {\bf 38} (1997), 1823-1842.
\vs.1 \nd
[Ke3]\ J.\ Kellendonk, The local structure of tilings and their integer
group of coinvariants, {\it Commun.\ Math.\ Phys.} {\bf 187} (1997), 115-157.
\vs.1 \nd
[Ke4]\ J.\ Kellendonk, On K$_0$-groups for substitution tilings, preprint.
\vs.1 \nd
[KeP]\ J.\ Kellendonk and I.\ Putnam, Tilings, C$^*$-algebras and K-theory,
preprint.
\vs.1 \nd
[PaT]\ W.\ Parry and S.\ Tuncel, {\it Classification Problems in Ergodic
Theory}, Lecture Note Series, {\bf 67}, Cambridge University Press, Cambridge, 1982.
\vs.1 \nd
[Put]\ I.\ Putnam, The ordered K-theory of C$^*$-algebras associated with
substitution tilings, preprint.
\vs.1 \nd
[Ra1]\ C.\ Radin, The pinwheel tilings of the plane, {\it Annals of
Math.}  {\bf 139} (1994), 661-702.
\vs.1 \nd
[Ra2]\ C.\ Radin, {\it Miles of Tiles}, Student Math. Lib. {\bf 1}, Amer.\ Math.\ Soc.,
Providence, 1999.
\vs.1 \nd
[RaS]\ C.\ Radin and L.\ Sadun, An algebraic invariant for substitution tiling systems, 
{\it Geometriae Dedicata} {\bf 73} (1998), 21-37. 
\vs1
\nd
[RaW]\ C.\ Radin and M.\ Wolff, Space tilings and local isomorphism,
{\it Geometriae Dedicata} {42} (1992), 355-360.
\vs.1 \nd
[Sol]\ B.\ Solomyak, Nonperiodicity implies unique composition for self-similar
translationally finite tilings, {\it Discrete Comput.\ Geom.} {\bf 20} (1998),
265-279.
\vs.1 \nd
[SuW]\ D.\ Sullivan and R.F.\ Williams, The homology of attractors, {\it Topology} {\bf 15}
(1976), 259-262. 
\vs1
\nd \line{Nicholas Ormes, Mathematics Department, University of Texas, Austin, TX\ \ 78712}

\nd {\it Email address}: {\tt ormes@math.utexas.edu}
\vs.1 \nd
\line{Charles Radin, Mathematics Department, University of Texas, Austin, TX\ \ 78712}

\nd {\it Email address}: {\tt radin@math.utexas.edu}
\vs.1 \nd
\line{Lorenzo Sadun, Mathematics Department, University of Texas, Austin, TX\ \ 78712}

\nd {\it Email address: {\tt sadun@math.utexas.edu}}
\vfill 
\end